\documentclass{amsart}
\usepackage{amssymb,amstext,amsmath,amscd,amsthm,amsfonts,enumerate,graphicx,latexsym,stmaryrd,multicol, mathrsfs, mathtools}

\usepackage[usenames]{color}
\usepackage[all]{xy}
\usepackage{extarrows}
\usepackage{mathabx}
\usepackage{cite}
\usepackage{textcomp}
\usepackage{hyperref}

\usepackage[top=25truemm,bottom=25truemm,left=25truemm,right=25truemm]{geometry}

\newtheorem{thm}{Theorem}[section]
\newtheorem{lem}[thm]{Lemma}
\newtheorem{prop}[thm]{Proposition}
\newtheorem{cor}[thm]{Corollary}
\theoremstyle{definition}
\newtheorem{dfn}[thm]{Definition}

\newtheorem{ques}[thm]{Question}

\newtheorem{ex}[thm]{Example}

\newtheorem{rmk}[thm]{Remark}
\theoremstyle{remark}

\newtheorem*{ac}{Acknowledgments}

\numberwithin{equation}{thm}
\def\add{\operatorname{add}}

\def\ass{\operatorname{Ass}}

\def\Cok{\operatorname{Coker}}

\def\Cone{\operatorname{Cone}}
\def\D{\mathrm{D}}
\def\dell{\operatorname{dell}}
\def\depth{\operatorname{depth}}
\def\E{\mathcal{E}}
\def\Ext{\operatorname{Ext}}
\def\End{\operatorname{End}}
\def\Endl{\operatorname{\underline{End}}}

\def\Gdim{\operatorname{G-dim}}
\def\ge{\geqslant}

\def\grade{\operatorname{grade}}
\def\grd{\mathsf{Grd}}

\def\Hom{\operatorname{Hom}}

\def\image{\operatorname{Im}}

\def\K{\mathcal{K}}
\def\Ker{\operatorname{Ker}}

\def\le{\leqslant}
\def\lhom{\operatorname{\underline{Hom}}}
\def\lmod{\operatorname{\underline{mod}}}
\def\m{\mathfrak{m}}

\def\mod{\operatorname{mod}}

\def\p{\mathfrak{p}}
\def\pd{\operatorname{pd}}

\def\proj{\operatorname{proj}}

\def\refl{\mathsf{Ref}}

\def\sgrade{\operatorname{s\text{.}grade}}
\def\Soc{\operatorname{Soc}}
\def\spec{\operatorname{Spec}}

\def\syz{\Omega}
\def\Syz{\operatorname{Syz}}
\def\t{\mathrm{T}}
\def\Trf{\mathsf{Trf}}

\def\tr{\operatorname{Tr}}
\def\trace{\operatorname{tr}}

\def\X{\mathcal{X}}

\def\Y{\mathcal{Y}}
\def\ZZ{\mathbb{Z}}
\begin{document}
\allowdisplaybreaks
\title[Auslander-type grade conditions and subcategories related to syzygy modules]{Auslander-type grade conditions and subcategories related to syzygy modules}
\author{Yuya Otake}
\address{Graduate School of Mathematical Sciences, The University of Tokyo, 3-8-1 Komaba, Meguro-ku, Tokyo 153-8914, Japan}
\email{yuyaotake31@g.ecc.u-tokyo.ac.jp}

\begin{abstract}
The question of when a syzygy module or a module satisfying a Serre condition satisfies the cohomological vanishing conditions known as higher torsionfreeness has been extensively studied. In this paper, we address these questions using grade conditions introduced and studied by Auslander--Bridger and Auslander--Reiten. We also investigate homological conditions on Ext modules and, as an application, give a partial answer over commutative rings to a question of Auslander--Bridger concerning whether higher reflexive modules are higher torsionfree.
\end{abstract}
\maketitle
\tableofcontents
\section{Introduction}

A fundamental problem in ring theory is to determine when a given module satisfies desirable homological conditions. In particular, understanding when higher syzygy modules satisfy certain cohomological vanishing conditions has been an important subject of study. Throughout this paper, all modules are assumed to be right modules unless otherwise specified.
Let $\Lambda$ be a two-sided noetherian ring. A finitely generated $\Lambda$-module $X$ is called {\em $n$-torsionfree} if $\Ext^k_{\Lambda^{\mathrm{op}}}(\tr X,\Lambda)=0$ for all $1\le k\le n$, where $\tr(-)$ denotes the Auslander transpose. This notion was introduced by Auslander and Bridger \cite{AB} as a natural higher version of torsionlessness. In particular, $1$-torsionfreeness and $2$-torsionfreeness coincide with torsionlessness and reflexivity, respectively. Higher torsionfree modules enjoy various remarkable properties arising from Auslander--Bridger theory. For example, an $n$-torsionfree module $X$ is $n$-reflexive, in the sense that the canonical unit morphism $X\to\syz^n\tr\syz^n\tr X$ is an isomorphism in the stable category, where $\syz^n X$ denotes an $n$-th syzygy of $X$ in a projective resolution. Higher torsionfree modules also play important roles in the representation theory of rings, particularly in the theory of Cohen--Macaulay representations; see \cite{LW,Yos90}.
Auslander and Bridger \cite{AB} investigated when a given $n$-syzygy module is $n$-torsionfree and discovered a close relationship between this question and Auslander's grade conditions on Ext modules. On the other hand, Serre conditions have long played a fundamental role in commutative algebra. If a commutative noetherian ring $R$ satisfies Serre's condition $(S_n)$, then every $n$-syzygy $R$-module also satisfies $(S_n)$. This naturally raises the converse question of when a module satisfying $(S_n)$ is an $n$-syzygy module or an $n$-torsionfree module. Such questions were classically investigated by Ischebeck \cite{Is}, Fossum and Reiten \cite{FR}, Foxby \cite{Fox}, and Evans and Griffith \cite{EG}. Further developments can be found in \cite{Masek,MTT}, among others.
The classical Serre condition compares the depth of a module with the height of a prime ideal. In our setting, however, it is more natural to compare the depth of a module with that of the ring itself. This leads us to the modified Serre condition, which can also be formulated for noetherian algebras over commutative noetherian rings. More precisely, let $\Lambda$ be a noetherian algebra over a commutative noetherian ring $R$. A finitely generated $\Lambda$-module $X$ is said to satisfy the {\em $n$-th modified Serre condition} over $\Lambda$ if, for any prime ideal $\p$ of $R$, the depth of $X_\p$ over $R_\p$ is at least the minimum of $n$ and the depth of $\Lambda_\p$ over $R_\p$.
Thus, we are led to consider the following inclusions of subcategories:
$$
\widetilde{S_n}(\Lambda)\supseteq\Syz_n(\Lambda)\supseteq\refl_n(\Lambda)\supseteq\Trf_n(\Lambda).
$$
Here, $\widetilde{S_n}(\Lambda)$ denotes the full subcategory of finitely generated $\Lambda$-modules satisfying the $n$-th modified Serre condition over $\Lambda$, while $\Syz_n(\Lambda)$, $\refl_n(\Lambda)$, and $\Trf_n(\Lambda)$ denote the subcategories of $n$-syzygy, $n$-reflexive, and $n$-torsionfree modules, respectively. Except for the leftmost inclusion, these inclusions hold over arbitrary two-sided noetherian rings.
Our main concern is to find conditions on individual modules under which the gaps between these subcategories can be bridged. One approach to this problem, originating in Auslander--Bridger theory and developed further by several authors, involves the finiteness of Gorenstein dimension at localizations; see, for example, \cite{AB,Masek,MJ}. In this paper, we approach these questions through grade conditions on Ext modules. Our first main result shows that such conditions provide characterizations of higher torsionfreeness for syzygy modules and modules satisfying modified Serre conditions.

\begin{thm}[Corollaries \ref{EGhighAB} and \ref{Serregrade}]\label{main1}
\begin{enumerate}[\rm(1)]
   \item Let $\Lambda$ be a two-sided noetherian ring and $n\ge0$ an integer. Then a finitely generated $\Lambda$-module $X$ is $n$-torsionfree if and only if $X$ is an $n$-th syzygy module and 
   $$
   \grade_\Lambda\Ext^k_{\Lambda^{\mathrm{op}}}(\tr X,\Lambda)\ge k-1\text{ for any }1\le k\le n.
   $$
   \item Let $\Lambda$ be a noetherian algebra over a commutative noetherian ring $R$, and let $n\ge0$ be an integer. Assume that $\Lambda$ satisfies the $n$-th modified Serre condition over $R$. Then a finitely generated $\Lambda$-module $X$ is $n$-torsionfree if and only if $X$ satisfies the $n$-th modified Serre condition over $\Lambda$ and
   $$
   \grade_R\Ext^k_{\Lambda^{\mathrm{op}}}(\tr X,\Lambda)\ge k \text{ for any }1\le k\le n.
   $$
\end{enumerate}
\end{thm}

Part (1) is essentially a consequence of Auslander--Bridger theory. It directly yields the remarkable left-right symmetry of Auslander--Reiten's condition $(d_n)$ established by Hoshino and Nishida \cite{HN} and Huang \cite{Hu08}; see Corollary \ref{symm}. It also gives a new description of the delooping level introduced by G\'elinas \cite{Gel}; see Proposition \ref{deloop} and Corollary \ref{deloop2}. Part (2) refines the classical result on higher torsionfreeness of maximal Cohen--Macaulay modules over orders with isolated singularities \cite{Aus78,EG,LW,Yos90}. In the commutative case, it also recovers a recent result of Mendoza-Rubio and Jorge-P\'erez \cite[Theorem 1.1]{MJ} on modules whose duals have locally finite Gorenstein dimension.

We next consider the extension-closedness of higher syzygy and torsionfree modules. This question has been extensively studied, notably by Auslander and Reiten \cite{AR1,AR2}, who related it to Auslander-type grade conditions; see also \cite{Huang99,HI,GT,DT}. Recently, Huang \cite{Hu24} investigated the relationship between these conditions and weakly Gorenstein modules studied by Ringel and Zhang \cite{RZ1,RZ2}, while Hanihara \cite{Han} related Auslander-type conditions to the quasi-abelian property of categories of reflexive modules.
We study extensions of individual higher torsionfree modules in more detail. Auslander and Reiten \cite{AR2} characterized, in terms of a strong grade condition, when every extension of a given $n$-torsionfree module $C$ by a projective module is again $n$-torsionfree. We show that, when the little finitistic dimension of $\Lambda$ is finite and $n$ is sufficiently large, even the intersection of the category of extensions of $C$ by projective modules with the category of $n$-torsionfree modules becomes ``trivial'', in the sense that such an extension is $n$-torsionfree only if it splits; see Theorem \ref{ARsplit}. This recovers several results of Dey and Takahashi \cite{DT}.

Finally, we consider the gap between $n$-reflexive modules and $n$-torsionfree modules, a question implicitly raised by Auslander and Bridger \cite{AB}. By definition, $1$-reflexivity coincides with torsionlessness. For $n=2$, Auslander and Bridger showed that $2$-reflexivity coincides with reflexivity over commutative noetherian rings. Very recently, Marczinzik \cite{Marc} showed that this is not the case for noncommutative algebras. Our second main result extends the commutative case to $n=3$.

\begin{thm}[Theorem \ref{3ref3trf}]\label{Main2}
Over a commutative noetherian ring, $3$-reflexivity and $3$-torsionfreeness coincide.
\end{thm}

To prove this theorem, we investigate the categories of modules isomorphic to $\Ext^n_{\Lambda^{\mathrm{op}}}(N,\Lambda)$ for some finitely generated $\Lambda^{\mathrm{op}}$-module $N$, and their submodules. Such categories were studied by Auslander in \cite{Aus69}, based on his theory of coherent functors \cite{AusC}. We first give a characterization of the subcategory consisting of submodules of such Ext modules. We also characterize $n$-reflexivity in terms of the projective dimension of a certain Ext module. We show that, when $\Lambda$ is commutative, every submodule of $\Ext^2_\Lambda(N,\Lambda)$ with projective dimension at most one is zero. This is a higher version of a classical result of Auslander \cite[Corollary 3.6]{Aus69} and leads to our second main theorem. Finally, we establish an Endo--Goto-type result for higher reflexivity over noetherian algebras over commutative noetherian rings; see \cite{EGref}.

The organization of this paper is as follows. In Section 2, we recall some fundamental notions from stable module theory and study the relationship between higher syzygies and torsionfreeness through grade conditions. As applications, we recover the left-right symmetry of Auslander--Reiten's condition $(d_n)$ and give a description of the delooping level. In Section 3, we investigate the relationship between modified Serre conditions and higher torsionfreeness over noetherian algebras. We also study how the local finiteness of G-dimension is related to grade conditions and recover some recent results on higher torsionfreeness. In Section 4, we study extensions of higher torsionfree modules and show that certain extensions must split when the torsionfreeness is sufficiently high. As applications, we recover and generalize some previous results, including noncommutative versions. Finally, in Section 5, we investigate submodules of Ext modules and give characterizations of the corresponding subcategories. As applications, we characterize higher reflexivity in terms of the projective dimensions of Ext modules and prove our second main theorem.
\section{Stable module theory and higher torsionfreeness}
Throughout this article, assume that all rings are two-sided noetherian, and that all modules are finitely generated right modules.
For a two-sided noetherian ring $\Lambda$, we denote by $\mod\Lambda$ the category of (finitely generated right) $\Lambda$-modules, and by $\proj\Lambda$ the full subcategory of $\mod\Lambda$ consisting of (finitely generated) projective $\Lambda$-modules.
We denote by $\lmod\Lambda$ the {\em stable category} of $\mod\Lambda$, that is, the ideal quotient of $\mod\Lambda$ by $\proj\Lambda$.
The objects of $\lmod\Lambda$ are the same as those of $\mod\Lambda$.  
For objects $X, Y$ in $\mod\Lambda$, the morphism set is defined by
\[
\Hom_{\lmod\Lambda}(X, Y) = \lhom_\Lambda(X, Y) = \Hom_\Lambda(X, Y)/\mathscr{P}(X, Y),
\]
where $\mathscr{P}(X, Y)$ denotes the subgroup of $\Hom_\Lambda(X, Y)$ consisting of $\Lambda$-homomorphisms that factor through some projective module.  
For any $\Lambda$-homomorphism $f : X \to Y$, we denote by $\underline{f}$ the image of $f$ in $\lhom_\Lambda(X, Y)$.
Let us recall the syzygy functor and the (Auslander) transpose functor, which play a fundamental role in the stable module theory.
We denote by $\Lambda^{\mathrm{op}}$ the opposite ring of $\Lambda$, and by $(-)^\ast$ the $\Lambda$-dual functor $\Hom_\Lambda(-, \Lambda)$.

\begin{dfn}
Let $M \in \mod\Lambda$, and take a projective presentation $P_1 \xrightarrow{\partial} P_0 \to M \to 0$ of $M$.
\begin{enumerate}[\rm(1)]
  \item
  We denote by $\syz M$ the image of $\partial$ and call it the {\em syzygy} of $M$.
  Then $\syz M$ is uniquely determined by $M$ up to projective summands.
  This assignment induces an additive functor $\syz : \lmod\Lambda \to \lmod\Lambda$.
  For each integer $n \ge 1$, the {\em $n$-th syzygy functor} $\syz^n$ is defined inductively by $\syz^n := \syz \circ \syz^{n-1}$.
  \item
  We denote by $\tr M$ the cokernel of the $\Lambda$-dual homomorphism $\partial^\ast$, and call it the {\em Auslander transpose} (or simply, the {\em transpose}) of $M$.
  Then $\tr M$ is uniquely determined by $M$ up to projective summands.
  This assignment induces an additive functor $\tr : \lmod\Lambda \to \lmod{\Lambda^{\mathrm{op}}}$, which is an anti-equivalence.
\end{enumerate}
\end{dfn}

Next, we prepare some notation concerning subcategories.

\begin{dfn}
Let $\X$ and $\Y$ be full subcategories of $\mod\Lambda$.
\begin{enumerate}[\rm(1)]
  \item
  We denote by $\add\X$ the {\em additive closure} of $\X$ in $\mod\Lambda$.
  Namely, $\add\X$ is the full subcategory of $\mod\Lambda$ consisting of all direct summands of finite direct sums of objects in $\X$.
  \item We denote by $\X\ast\Y$ the subcategory of $\mod\Lambda$ consisting of modules $Z$ such that there exists an exact sequence $0\to X\to Z\to Y\to 0$ in $\mod\Lambda$ with $X\in \X$ and $Y\in\Y$.
  \item
  Let $n \ge 0$ be an integer.
  We denote by $\Syz_n\X$ the full subcategory of $\mod\Lambda$ consisting of $n$-th syzygies of objects in $\X$; that is, $\Lambda$-modules $M$ for which there exists an exact sequence
  $0 \to M \to P_{n-1} \to \cdots \to P_0 \to X \to 0$ of $\Lambda$-modules such that $P_i \in \proj\Lambda$ for all $0 \le i \le n-1$, and $X \in \X$.
  Also, we denote $\Syz_n(\mod\Lambda)$ simply by $\Syz_n(\Lambda)$.
  Similarly, we denote by $\tr\X$ the full subcategory of $\mod{\Lambda^{\mathrm{op}}}$ consisting of transposes of objects in $\X$; that is, $\Lambda^{\mathrm{op}}$-modules $N$ for which there exists an exact sequence $P_1 \xrightarrow{\partial} P_0 \to X \to 0$ of $\Lambda$-modules such that $P_0, P_1 \in \proj\Lambda$, $X \in \X$, and the cokernel of the $\Lambda$-dual homomorphism $\partial^\ast$ is isomorphic to $N$.
\end{enumerate}
\end{dfn}

Let us recall a well-known and fundamental fact in stable module theory.  
Let $n \ge 0$ be an integer.  
It is explicitly stated in \cite[Corollary 3.4]{AR2} that the pair of functors $(\tr \syz^n \tr, \syz^n)$ on $\lmod\Lambda$ forms an adjoint pair.  
That is, there exist functorial isomorphisms
$$
\theta^n_{M, N} : \lhom_{\Lambda}(M, \syz^n N) \xrightarrow{\sim} \lhom_\Lambda(\tr \syz^n \tr M, N)
$$
for all $M, N \in \lmod\Lambda$; see also \cite{AB, GT, Omorp}.
We denote the counit of this adjoint pair by $\psi^n_{(-)} : \tr \syz^n \tr \syz^n(-) \to (-)$, and the unit by $\varphi^n_{(-)} : (-) \to \syz^n \tr \syz^n \tr(-)$.
That is, for each $M \in \lmod\Lambda$, the morphisms $\underline{\psi^n_M}$ and $\underline{\varphi^n_M}$ in $\lmod\Lambda$ are given by
$$
\underline{\psi^n_M} = \theta^n_{\syz^n M, M}(\underline{1_{\syz^n M}}) : \tr \syz^n \tr \syz^n M \to M
$$
and
$$
\underline{\varphi^n_M} = \left(\theta^n_{M, \tr \syz^n \tr M}\right)^{-1}(\underline{1_{\tr \syz^n \tr M}}) : M \to \syz^n \tr \syz^n \tr M.
$$

The notions of $n$-torsionfreeness, $n$-reflexivity, and Rees's grade play a central role in the Auslander--Bridger theory \cite{AB}.

\begin{dfn}
Let $M \in \mod\Lambda$ and let $n \ge 0$ be an integer.
\begin{enumerate}[\rm(1)]
  \item
  The module $M$ is said to be {\em $n$-torsionfree} if $\Ext^i_{\Lambda^{\mathrm{op}}}(\tr M, \Lambda) = 0$ for all $1 \le i \le n$.
  We denote by $\Trf_n(\Lambda)$ the category of $n$-torsionfree modules.
  \item
  The module $M$ is said to be {\em $n$-reflexive} if the unit morphism $\underline{\varphi^n_M} : M \to \syz^n \tr \syz^n \tr M$ is an isomorphism in the stable module category $\lmod\Lambda$.
  We denote by $\refl_n(\Lambda)$ the category of $n$-reflexive modules.
  \item
  The {\em grade} of $M$, written $\grade_\Lambda M$, is defined to be the smallest integer $i \ge 0$ such that $\Ext^i_{\Lambda}(M, \Lambda) \ne 0$.
  We denote by $\grd_n(\Lambda)$ the category of modules whose grade is greater than or equal to $n$.
\end{enumerate}
\end{dfn}

Let us see that $1$-torsionfreeness and $2$-torsionfreeness are equivalent to the well-known notions of torsionlessness and reflexivity, respectively.
Let $M$ be a $\Lambda$-module.  
We consider the natural $\Lambda$-homomorphism $\varphi_M : M \to M^{\ast\ast}$ defined by $\varphi_M(x)(f) = f(x)$ for $x \in M$ and $f \in M^\ast$.  
The module $M$ is called {\em torsionless} if the homomorphism $\varphi_M$ is injective, and {\em reflexive} if it is bijective.
For this homomorphism $\varphi_M$, there exists a short exact sequence called the {\em Auslander sequence}:
$$
0 \to \Ext^1_{\Lambda^{\mathrm{op}}}(\tr M, \Lambda) \to M \xrightarrow{\varphi_M} M^{\ast\ast} \to \Ext^2_{\Lambda^{\mathrm{op}}}(\tr M, \Lambda) \to 0,
$$
see \cite{AusC, AB}.  
Thus, $M$ is torsionless if and only if it is $1$-torsionfree, and reflexive if and only if it is $2$-torsionfree.
As stated in \cite[Appendix]{AB}, the image $\image\varphi_M$ is naturally isomorphic to $\syz\tr\syz\tr M$ in $\lmod\Lambda$, and the natural homomorphism $M \to \image\varphi_M$ induced by $\varphi_M$ coincides with $\underline{\varphi^1_M} : M \to \syz\tr\syz\tr M$ in $\lmod\Lambda$.
Similarly, since $(-)^\ast$ is naturally isomorphic to $\syz^2 \tr(-)$, the morphism $\underline{\varphi_M} : M \to M^{\ast\ast}$ coincides with the unit morphism $\underline{\varphi^2_M} : M \to \syz^2 \tr \syz^2 \tr M$ in $\lmod\Lambda$.
For any integer $n\ge0$, the following duality for $n$-torsionfree modules is known in Auslander--Bridger theory \cite[Theorem 2.17]{AB}.

\begin{thm}[Auslander--Bridger]\label{trfdual}
Let $n\ge0$ be an integer. Then there exists a duality between stable module categories given by
$$
\xymatrix{
\underline{\Trf_n}(\Lambda)\ar@<.5mm>[rr]^-{\syz^n\tr}
&&\underline{\Trf_n}(\Lambda^{\mathrm{op}})\ar@<.5mm>[ll]^-{\syz^n\tr}
}
$$
In particular, if $M$ is an $n$-torsionfree $\Lambda$-module, then the $k$-th unit morphism $\underline{\varphi^k_M}:M\to\syz^k\tr\syz^k\tr M$ is an isomorphism in $\lmod\Lambda$ for any $1\le k\le n$. Namely, $M$ is $k$-reflexive for any $1\le k\le n$.
\end{thm}

Since $n$-reflexive modules are clearly $n$-syzygies, for any $n \ge 0$, we have inclusions of subcategories
$$
\add\Syz_n(\Lambda) \supseteq \Syz_n(\Lambda) \supseteq \refl_n(\Lambda) \supseteq \Trf_n(\Lambda)
$$
of $\mod\Lambda$. It is an interesting question when each of these inclusions becomes an equality, and this question has been extensively studied. 
As is well known, a module is torsionless if and only if it is a first syzygy.  
Therefore, any direct summand of a first syzygy is again a first syzygy.
Hence, the equalities $\add\Syz_1(\Lambda) = \Syz_1(\Lambda) = \refl_1(\Lambda) = \Trf_1(\Lambda)$ hold. 
In higher cases, direct summands of higher syzygies can be characterized by the splitting of the unit morphism $\underline{\varphi^n_{(-)}}$ as follows.

\begin{lem}\label{EGaddsyz}
Let $M$ be a $\Lambda$-module and let $n \ge 0$ be an integer.  
Then $M$ belongs to $\add\Syz_n(\Lambda)$ if and only if the unit morphism $\underline{\varphi^n_M} : M \to \syz^n \tr \syz^n \tr M$ is a split monomorphism in $\lmod\Lambda$.
Furthermore, in the case where $n = 2$, $M$ belongs to $\add\Syz_2(\Lambda)$ if and only if the natural homomorphism $\varphi_M : M \to M^{\ast\ast}$ is a split injection in $\mod\Lambda$.
\end{lem}

\begin{proof}
Since the “if” part is clear, it suffices to prove the “only if” part.  
We handle here the case where $n = 2$.  
Since $M$ is an object of $\add\Syz_2(\Lambda)$, there exists a $\Lambda^{\mathrm{op}}$-module $N$ such that $M$ is a direct summand of the $\Lambda$-dual $N^\ast$.  
Namely, we can take a split injection $f : M \to N^\ast$ in $\mod\Lambda$.  
Also, the diagram
$$
\xymatrix@R-1pc@C-1pc{
M\ar[rrr]^-{f}\ar[dd]^-{\varphi_M}&&&N^\ast\ar[dd]^-{\varphi_{N^\ast}}\\
&&&\\
M^{\ast\ast}\ar[rrr]^-{f^{\ast\ast}}&&&N^{\ast\ast\ast}
}
$$
is commutative in $\mod\Lambda$.  
Note that $\varphi_{N^\ast}$ is a split injection since the composite $(\varphi_N)^\ast \circ \varphi_{N^\ast}$ equals the identity map on $N^\ast$ in $\mod\Lambda$.  
As both $f$ and $\varphi_{N^\ast}$ are split injections in $\mod\Lambda$, so is $\varphi_M$.
The remaining assertion in the case where $n \ne 2$ can be proved in essentially the same way.  
Assume that $M$ belongs to $\add\Syz_n(\Lambda)$ for $n\ge1$.  
This assumption implies that there exist a $\Lambda$-module $U$ and a split monomorphism $\underline{g} : M \to \syz^n U$ in $\lmod\Lambda$.  
Then we obtain the following commutative diagram in $\lmod\Lambda$:
$$
\xymatrix@R-1pc@C-1pc{
M \ar[rrrrr]^-{\underline{g}} \ar[dd]^-{\underline{\varphi^n_M}} &&&&& \syz^n U \ar[dd]^-{\underline{\varphi^n_{\syz^n U}}} \\
&&&&& \\
\syz^n \tr \syz^n \tr M \ar[rrrrr]^-{\underline{\syz^n \tr \syz^n \tr g}} &&&&& \syz^n \tr \syz^n \tr \syz^n U.
}
$$
As shown in \cite[Page 63]{AB}, the composite $\underline{\syz^n \psi^n_{U}} \circ \underline{\varphi^n_{\syz^n U}}$ is equal to the identity morphism on $\syz^n U$ in $\lmod\Lambda$.  
Hence, $\underline{\varphi^n_{\syz^n U}}$ is a split monomorphism in $\lmod\Lambda$.  
Since both $\underline{g}$ and $\underline{\varphi^n_{\syz^n U}}$ are split monomorphisms in $\lmod\Lambda$, it follows that $\underline{\varphi^n_M}$ is also a split monomorphism.
\end{proof}

It is a classical and well-studied problem to determine when $n$-th syzygy modules are $n$-torsionfree.
Auslander and Bridger related this question to grade theory.
The following theorem, which is precisely \cite[Proposition 2.26]{AB}, is a fundamental result in this direction.

\begin{thm}[Auslander--Bridger]\label{EGABclass}
Let $n > 0$ be an integer and $M\in\mod\Lambda$.  
Then the following conditions are equivalent.
\begin{enumerate}
   \item The $(k+1)$-st syzygy $\syz^{k+1} M$ is $(k+1)$-torsionfree for any $1\le k\le n$.
   \item One has $\grade_{\Lambda^{\mathrm{op}}}\Ext^{k+1}_\Lambda(M,\Lambda)\ge k$ for any $1\le k\le n$.
\end{enumerate}
\end{thm}

Auslander and Reiten \cite{AR1} called $\Lambda$ a ring satisfying condition $(d_n)$ if every $M\in\mod\Lambda$ satisfies condition (2) of the above theorem. This condition is denoted by $g_n(1)$ in Huang--Iyama \cite{HI}. By the above theorem, this condition is equivalent to the equality $\Syz_k(\Lambda)=\Trf_k(\Lambda)$ for any $1\le k\le n+1$, which is also equivalent to $\add\Syz_k(\Lambda)=\Syz_k(\Lambda)=\Trf_k(\Lambda)$ for any $1\le k\le n+1$.
Furthermore, Auslander and Bridger in \cite[Proposition 4.21]{AB} proved that when $\Lambda$ is commutative, condition $(d_n)$ is also equivalent to the equality $\syz^{n+1}(\mod\Lambda)=\Trf_{n+1}(\Lambda)$, which holds if and only if $\Lambda$ is locally Gorenstein at all prime ideals of depth less than $n$.


We consider when a given direct summand of an $n$-th syzygy module is $n$-torsionfree. By combining a characterization of such direct summands with the higher transpose duality, we obtain the following.

\begin{lem}\label{highABlem}
Let $n\ge0$ be an integer and $M\in\add\Syz_n(\Lambda)$. Then $M$ is an $n$-torsionfree $\Lambda$-module if and only if $\syz^n\tr M$ is an $n$-torsionfree $\Lambda^{\mathrm{op}}$-module.
\end{lem}

\begin{proof}
By Theorem \ref{trfdual}, the ``only if'' part is true. We prove the converse. Since $M\in\add\Syz_n(\Lambda)$, Lemma \ref{EGaddsyz} shows that the unit $\underline{\varphi^n_M}:M\to\syz^n\tr\syz^n\tr M$ is a split monomorphism in $\lmod\Lambda$. Since $\syz^n\tr M$ is $n$-torsionfree, Theorem \ref{trfdual} again implies that $\syz^n\tr\syz^n\tr M$ is also $n$-torsionfree. Hence, $M$, being a direct summand of this module in $\lmod\Lambda$, is also $n$-torsionfree.
\end{proof}

As a consequence, we obtain the following result, which shows that the gap between $n$-torsionfreeness and $n$-syzygies can be bridged by an Auslander--Bridger-type grade condition on the transpose.

\begin{cor}\label{EGhighAB}
Let $n\ge0$ be an integer and $M\in\mod\Lambda$.
Then $M$ is $n$-torsionfree if and only if $M\in\add\Syz_n(\Lambda)$ and the following condition is satisfied:
$$
\grade_{\Lambda}\Ext^k_{\Lambda^{\mathrm{op}}}(\tr M,\Lambda)\ge k-1\text{ for all }1\le k\le n.
$$
\end{cor}

\begin{proof}
The ``only if'' part is trivial. We prove the ``if'' part. By the grade condition on the transpose, Theorem \ref{EGABclass} implies, in particular, that $\syz^n\tr M$ is $n$-torsionfree. Hence, the assertion follows from Lemma \ref{highABlem}.
\end{proof}

This result immediately recovers the remarkable left-right symmetry of condition $(d_n)$ established by Hoshino and Nishida \cite[Theorem 4.6]{HN} and Huang \cite[Theorem 2.4]{Hu08}.

\begin{cor}[Hoshino--Nishida, Huang]\label{symm}
Let $\Lambda$ be a two-sided noetherian ring and $n\ge0$ an integer. Then $\Lambda$ satisfies Auslander--Reiten's condition $(d_n)$ if and only if so does $\Lambda^{\mathrm{op}}$.
\end{cor}

\begin{proof}
Suppose that $\Lambda$ satisfies condition $(d_n)$. By Theorem \ref{EGABclass}, it is enough to show that $\Syz_k(\Lambda^{\mathrm{op}})=\Trf_k(\Lambda^{\mathrm{op}})$ for any $1\le k\le n+1$. However, since $\Lambda$ satisfies condition $(d_n)$, that is, condition (2) of Theorem \ref{EGABclass}, every object in $\Syz_k(\Lambda^{\mathrm{op}})$ satisfies the grade condition on its transpose required in Corollary \ref{EGhighAB}, and hence is $k$-torsionfree.
\end{proof}

We now give an application to the delooping level introduced by Gélinas \cite{Gel}. For a module $M\in\mod\Lambda$, the {\em delooping level} of $M$, denoted by $\dell M$, is defined to be the infimum of integers $n\ge0$ such that $\syz^n M$ is a direct summand of $\syz^{n+1}N$ in $\lmod\Lambda$ for some $N\in\mod\Lambda$. This invariant is closely related to the existence of $n$-FPD hulls, introduced in \cite{OExt, OMart} as a refinement of FPD hulls over Gorenstein rings studied in \cite{ABu}. An {\em $n$-FPD hull} of $M$ is an exact sequence $0\to M\to Y\to X\to0$ such that $\pd_\Lambda Y\le n$ and $\Ext^i_\Lambda(X,\Lambda)=0$ for all $1\le i\le n+1$. The following result refines \cite[Theorem 1.3]{OExt}, where $\Lambda$ was assumed to be commutative and $M$ to have locally finite G-dimension. By Corollary \ref{EGhighAB}, these assumptions can be replaced by a grade condition.

\begin{prop}\label{deloop}
Let $\Lambda$ be a two-sided noetherian ring and $M\in\mod\Lambda$. Then
$$
\dell M\le\inf\{n\ge0\mid M\text{ admits an }n\text{-FPD hull }\}.
$$
Moreover, suppose that $d=\dell M<\infty$ and set $M_d=\syz^d M$.
Then equality holds if
$$
\grade_\Lambda\Ext^k_{\Lambda^{\mathrm{op}}}(\tr M_d,\Lambda)\ge k-1\text{ for any }1\le k\le d+1.
$$
\end{prop}

\begin{proof}
The first assertion was proved in \cite[Theorem 2.3]{OExt} under the assumption that $\Lambda$ is commutative. As noted in \cite[Remark 3.11(2)]{OMart}, this assumption is unnecessary. We prove the second assertion. Since $d=\dell M$, we have $M_d\in\add\Syz_{d+1}(\Lambda)$. Hence, Corollary \ref{EGhighAB} shows that $M_d$ is $(d+1)$-torsionfree. Thus, $M$ admits a $d$-FPD hull, and the assertion follows.
\end{proof}
We say that $\Lambda$ satisfies $(d_\infty)$ if it satisfies $(d_n)$ for all $n\ge0$.
\begin{cor}\label{deloop2}
Let $\Lambda$ be a two-sided noetherian ring satisfying condition $(d_\infty)$. Then, for any $M\in\mod\Lambda$,
$$
\dell M=\inf\{n\ge0\mid M\text{ admits an }n\text{-FPD hull }\}.
$$
\end{cor}

\section{On modified Serre conditions}
Throughout this section, let $R$ be a commutative noetherian ring and let $\Lambda$ be a noetherian $R$-algebra, that is, a module-finite $R$-algebra, not necessarily commutative.
In this section, we study the relationship between Serre conditions and $n$-torsionfreeness. The classical Serre conditions relate the depth of a module to the heights of prime ideals.
However, in our homological setting, it is natural to consider the following modified version, which coincides with the classical Serre condition in the Cohen--Macaulay case.
We recall the definition of depth over a commutative local ring.
Let $R$ be a commutative local ring with residue field $k$ and $M$ an $R$-module.
The {\em depth} of $M$, denoted by $\depth_R M$, is defined by $\depth_R M=\inf\{i\ge0\mid\Ext^i_R(k,M)\ne0\}$.

\begin{dfn}\label{modSer}
Let $\Lambda$ be a noetherian $R$-algebra and $M$ a finitely generated $\Lambda$-module. We say that $M$ satisfies the \textit{modified Serre condition} $(\widetilde{S_n^R})$ if $\depth_{R_\p} M_\p \ge \min\{n,\depth_{R_\p} \Lambda_\p\}$ for every $\p\in\spec R$.
We denote by $\widetilde{S_n^R}(\Lambda)$ the full subcategory of $\mod\Lambda$ consisting of modules satisfying the modified Serre condition $(\widetilde{S_n^R})$.
\end{dfn}

In the case where $\Lambda=R$, this condition coincides with Foxby's condition $(b_n)$ \cite{Fox}. In later work, including that of Ma\c{s}ek \cite{Masek}, this condition is denoted by $(\widetilde{S}_n)$. We use the same notation in this paper.

\begin{lem}\label{syzsn}
Suppose that a finitely generated $\Lambda$-module $M$ is an $n$-th syzygy over $\Lambda$. Then $M$ satisfies the modified Serre condition $(\widetilde{S_n^R})$. In other words, the inclusion $\Syz_n(\Lambda)\subseteq \widetilde{S_n^R}(\Lambda)$ holds.
\end{lem}
\begin{proof}
Fix $\p\in\spec R$ and replace $R$ by $R_\p$. Since every $P\in\proj\Lambda$ is a direct summand of a finite direct sum of copies of $\Lambda$, one has $\depth_R P\ge \depth_R\Lambda$. Hence the assertion follows immediately from the depth lemma.
\end{proof}

The preceding lemma motivates us to ask when the modified Serre condition implies higher torsionfreeness. The following theorem gives such a criterion in terms of the associated primes of certain Ext modules.

\begin{thm}\label{SerreAss}
Let $n\ge0$ be an integer.
Let $M$ be a finitely generated $\Lambda$-module satisfying the modified Serre condition $(\widetilde{S_n^R})$. Suppose that, for every integer $1\le k\le n$, one has $$\ass_R\Ext^k_{\Lambda^{\mathrm{op}}}(\tr M,\Lambda)\subseteq\{\p\in\spec R\mid \depth_{R_\p}\Lambda_\p\ge k\}.$$
Then $M$ is an $n$-torsionfree $\Lambda$-module.
\end{thm}

\begin{proof}
First, we consider the case $n=1$. We show that $\Ext^1_{\Lambda^{\mathrm{op}}}(\tr M,\Lambda)=0$.
If not, the set of associated primes of $\Ext^1_{\Lambda^{\mathrm{op}}}(\tr M,\Lambda)$ is nonempty.
Take $\p\in\ass_R\Ext^1_{\Lambda^{\mathrm{op}}}(\tr M,\Lambda)$.
By the Auslander--Bridger sequence, $\Ext^1_{\Lambda^{\mathrm{op}}}(\tr M,\Lambda)$ is a submodule of $M$ as a $\Lambda$-module, and hence also as an $R$-module.
Thus $\p\in\ass_R M$, and so $\depth_{R_\p}M_\p=0$. Since $M$ satisfies $(\widetilde{S_1^R})$, we have $0=\depth_{R_\p}M_\p\ge\min\{1,\depth_{R_\p}\Lambda_\p\}$. Hence $\depth_{R_\p}\Lambda_\p=0$. On the other hand, since $\p$ is an associated prime of $\Ext^1_{\Lambda^{\mathrm{op}}}(\tr M,\Lambda)$, the hypothesis gives $\depth_{R_\p}\Lambda_\p\ge1$.
This is a contradiction.
Next, we consider the case $n=2$. By the case $n=1$, the module $M$ is torsionless.
We show that $\Ext^2_{\Lambda^{\mathrm{op}}}(\tr M,\Lambda)=0$.
If not, take $\p$ in the set of associated primes of $\Ext^2_{\Lambda^{\mathrm{op}}}(\tr M,\Lambda)$.
By the hypothesis, we have $\depth_{R_\p}\Lambda_\p\ge2$. Since $M$ satisfies the condition $(\widetilde{S_2^R})$, we also have $\depth_{R_\p}M_\p\ge2$.
On the other hand, the Auslander--Bridger sequence gives an exact sequence $0\to M\to M^{\ast\ast}\to\Ext^2_{\Lambda^{\mathrm{op}}}(\tr M,\Lambda)\to0$.
Since $M^{\ast\ast}$ is a second syzygy module, and hence in particular a first syzygy module, Lemma~\ref{syzsn} gives $\depth_{R_\p}(M^{\ast\ast})_\p\ge1$. Since $\p$ is an associated prime of $\Ext^2_{\Lambda^{\mathrm{op}}}(\tr M,\Lambda)$, its localization at $\p$ has depth zero. Hence the depth lemma gives $\depth_{R_\p}M_\p=1$. This is a contradiction. So we can conclude that $M$ is reflexive.
The cases $n=1$ and $n=2$ have been proved. We now consider the case $n\ge3$. By induction, we may assume that $M$ is $(n-1)$-torsionfree, and hence reflexive.
Take a projective resolution
$$\cdots\to Q_{n-2}\to Q_{n-3}\to\cdots\to Q_1\to Q_0\to M^\ast\to0$$
of the $\Lambda^{\mathrm{op}}$-module $M^\ast$.
Since $M$ is $(n-1)$-torsionfree, taking the $\Lambda$-dual and using $M\cong M^{\ast\ast}$, we obtain an exact sequence
$$0\to M\to Q_0^\ast\to Q_1^\ast\to\cdots\to Q_{n-3}^\ast\to Q_{n-2}^\ast\to\tr\syz^{n-3}(M^\ast)\to0.$$
Note that $\tr\syz^{n-3}(M^\ast)\cong\tr\syz^{n-1}\tr M$.
We claim that $\Ext^n_{\Lambda^{\mathrm{op}}}(\tr M,\Lambda)=0$. Assume that this module is nonzero, and choose $\p$ from its set of associated primes.
We have $\Ext^n_{\Lambda^{\mathrm{op}}}(\tr M,\Lambda)\cong\Ext^1_{\Lambda^{\mathrm{op}}}(\tr\tr\syz^{n-1}\tr M,\Lambda)\hookrightarrow\tr\syz^{n-1}\tr M\cong\tr\syz^{n-3}(M^\ast)$.
Hence $\p\in\ass_R\tr\syz^{n-3}(M^\ast)$.
By the hypothesis, $\depth_{R_\p}\Lambda_\p\ge n$. Since $M$ satisfies $(\widetilde{S_n^R})$, we also have $\depth_{R_\p}M_\p\ge n$.
Moreover, each $Q_i^\ast$ is a projective $\Lambda$-module, and hence $\depth_{R_\p}(Q_i^\ast)_\p\ge n$.
On the other hand, $\depth_{R_\p}(\tr\syz^{n-3}(M^\ast))_\p=0$.
Repeated use of the depth lemma gives $\depth_{R_\p}M_\p=n-1$, contrary to $\depth_{R_\p}M_\p\ge n$.
Thus $\Ext^n_{\Lambda^{\mathrm{op}}}(\tr M,\Lambda)=0$, and hence $M$ is $n$-torsionfree.
\end{proof}

The following corollary follows immediately from this theorem. It shows that Serre conditions and higher torsionfreeness are also related through Auslander-type grade conditions.

\begin{cor}\label{Serregrade}
Let $M$ be a finitely generated $\Lambda$-module and $n\ge0$ an integer.
Assume that $M$ satisfies the modified Serre condition $(\widetilde{S_n^R})$ over $\Lambda$ and that $\Lambda$
satisfies the modified Serre condition $(\widetilde{S_n^R})$ over $R$; that is, $M\in\widetilde{S_n^R}(\Lambda)$ and
$\Lambda\in\widetilde{S_n^R}(R)$.
If $$\grade_R\Ext^k_{\Lambda^{\mathrm{op}}}(\tr M,\Lambda)\ge k \text{ for every }1\le k\le n,$$ then $M$ is an $n$-torsionfree $\Lambda$-module.
\end{cor}

\begin{proof}
Take any integer $1\le k\le n$ and a prime ideal $\p\in\spec R$ such that $\depth_{R_\p}\Lambda_\p<k$.
Since $\Lambda$ satisfies the modified Serre condition over $R$, we have $\depth R_\p<k$.
Hence the grade assumption implies that $\Ext^k_{\Lambda^{\mathrm{op}}}(\tr M,\Lambda)_\p=0$.
Thus the desired depth condition is automatically satisfied, and the assertion follows from Theorem \ref{SerreAss}.
\end{proof}

We consider when the condition $\Lambda\in\widetilde{S_n^R}(R)$ is satisfied. Of course, when $\Lambda=R$, this condition is automatically satisfied.
More generally, if $\Lambda$ is a projective $R$-module, then the condition is also satisfied.
Another typical example where the above condition is satisfied is when $R$ is a commutative Cohen--Macaulay ring
and $\Lambda$ is an \textit{$R$-order}, that is, $\Lambda_\p$ is a maximal Cohen--Macaulay $R_\p$-module for every $\p\in\spec R$.

\begin{rmk}
\begin{enumerate}[\rm(1)]
   \item We consider the case where $\Lambda=R$ in Corollary \ref{Serregrade}, and the grade condition $\grade_R\Ext^k_R(\tr M,R)\ge k$ for every $1\le k\le n$.
   In \cite[Theorem 6.13(1)(2)]{EGGHIKMT}, Endo et al.\ prove a similar result by replacing this condition with the condition that $M_\p$ is $k$-torsionfree for every $1\le k\le n$ and every $\p\in\spec R$ with $\depth R_\p<k$.
   Their theorem is a higher version of \cite[Proposition 1.4.1]{BH}.
   Their local $k$-torsionfreeness condition is stronger than our grade condition, and the two conditions coincide when $n=1$.
   However, they are not necessarily equivalent when $n\ge2$.
   For example, when $n=2$, the difference can be easily seen.
   Let $R$ be a regular local ring of dimension $d\ge1$, and let $x$ be part of a regular system of parameters of $R$.
   Set $M=R/(x)$. Then $\grade_R\Ext^1_R(\tr M,R)=\grade_R R/(x)=1$ and $\Ext^2_R(\tr M,R)=0$.
   Thus the grade condition for $n=2$ is satisfied. However, the prime ideal $\p=(x)$ satisfies $\depth R_\p=1$, and $M_\p\cong R_\p/\p R_\p$ is not a torsionless $R_\p$-module.
   Hence the local reflexivity condition is not satisfied.
   \item Corollary \ref{Serregrade} can be regarded as a refinement of the following well-known fact in Cohen--Macaulay representation theory:
   Let $R$ be a $d$-dimensional commutative Cohen--Macaulay local ring and $\Lambda$ an $R$-order.
   If $M$ is a maximal Cohen--Macaulay $\Lambda$-module which is locally free on the punctured spectrum of $R$, then $M$ is a $d$-torsionfree $\Lambda$-module; see \cite{Aus78,EG,LW,Yos90}.
   Indeed, in this case, $M$ satisfies the modified Serre condition $(\widetilde{S_d^R})$, and $\Ext^k_{\Lambda^{\mathrm{op}}}(\tr M,\Lambda)$ has grade at least $d$ over $R$ for every $1\le k\le d$. See also \cite[Proposition 4.3]{DKLO}.
\end{enumerate}
\end{rmk}

Next, we consider the relationship between the local finiteness of G-dimension and Auslander-type grade conditions.
The notion of G-dimension was introduced by Auslander and Bridger \cite{AB} as a Gorenstein analogue of projective dimension.
Let $\Lambda$ be a two-sided noetherian ring.
A finitely generated $\Lambda$-module $X$ is called \textit{totally reflexive} if $\Ext^i_\Lambda(X,\Lambda)=0$ and $\Ext^i_{\Lambda^{\mathrm{op}}}(\tr X,\Lambda)=0$ for every $i\ge1$.
For a finitely generated $\Lambda$-module $M$, the \textit{G-dimension} of $M$, denoted by $\Gdim_\Lambda M$, is defined to be the infimum of all integers $n\ge0$ for which there exists an exact sequence
$
0\to X_n\to\cdots\to X_1\to X_0\to M\to0
$
with each $X_i$ totally reflexive.
If no such sequence exists, we set $\Gdim_\Lambda M=\infty$.
By \cite[Theorem 3.13]{AB}, it is known that $\Gdim_\Lambda M\le n$ if and only if $\syz^n M$ is totally reflexive.

As we show below, when $\Lambda=R$ is a commutative noetherian ring, the local finiteness of G-dimension implies grade conditions on Ext modules.
The following proposition is a refinement of \cite[Proposition 4.16 and Corollary 4.17]{AB}.

\begin{prop}\label{Gdimgrd}
Let $n\ge0$ be an integer, and assume that a finitely generated $R$-module $M$ satisfies the condition that $\Gdim_{R_\p}M_\p<\infty$ for every $\p\in\spec R$ with $\depth R_\p<n$. Then $\grade_R\Ext^m_R(M,R)\ge k$ for any $1\le k\le n$ and $m\ge k$.
\end{prop}

\begin{proof}
Fix integers $1\le k\le n$ and $m\ge k$ and take a prime ideal $\p$ of $R$ with $\depth R_\p<k$. It is enough to show that $\Ext^m_R(M,R)_\p=0$.
By the assumption and the Auslander--Bridger formula, $M_\p$ has finite G-dimension over $R_\p$ at most $\depth R_\p$, which is less than $k$.
Hence $\Ext^m_R(M,R)_\p\cong\Ext^m_{R_\p}(M_\p,R_\p)=0$. This completes the proof.
\end{proof}

Combining the main result of this section with the above proposition, we obtain the following corollaries. They show that higher torsionfreeness follows from being a higher syzygy or satisfying the modified Serre condition, under the assumption that the dual module has locally finite G-dimension. The second statement is precisely a recent result of Mendoza-Rubio and Jorge-P\'erez \cite[Theorem 1.1]{MJ}.

\begin{cor}\label{MJtype}
Let $M$ be a finitely generated $R$-module and $n\ge0$ an integer.
\begin{enumerate}[\rm(1)]
\item Assume that $\Gdim_{R_\p}(M_\p)^*<\infty$ for every $\p\in\spec R$ with $\depth R_\p<n$. Then $M\in\add\Syz_{n+1}(R)$ if and only if $M\in\Trf_{n+1}(R)$.
\item Assume that $\Gdim_{R_\p}(M_\p)^*<\infty$ for every $\p\in\spec R$ with $\depth R_\p<n$. Then $M\in\widetilde{S_n}(R)$ if and only if $M\in\Trf_n(R)$.
\end{enumerate}
\end{cor}

\begin{proof}
First, note that the assumption is equivalent to the local finiteness of the G-dimension of $\tr M$, since $\syz^2\tr M\cong M^*$ up to projective summands. In both assertions, the ``only if'' direction is the only one that needs to be proved.
\begin{enumerate}[\rm(1)]
\item Assume that $M\in\add\Syz_{n+1}(R)$. By the assumption and Proposition \ref{Gdimgrd}, we have $\grade_R\Ext^{k+1}_R(\tr M,R)\ge k$ for every $1\le k\le n$. Hence Corollary \ref{EGhighAB} implies that $M$ is $(n+1)$-torsionfree.
\item Assume that $M\in\widetilde{S_n}(R)$. By the assumption and the preceding lemma, we have $\grade_R\Ext^k_R(\tr M,R)\ge k$ for every $1\le k\le n$. Hence Corollary \ref{Serregrade} implies that $M$ is $n$-torsionfree.
\end{enumerate}
\end{proof}

In the following remark, we also consider the finiteness of G-dimension of $M$ itself.

\begin{rmk}\label{nondual}
In each of (1) and (2) of Corollary \ref{MJtype}, the assertion remains true if the assumption ``$\Gdim_{R_\p}(M_\p)^*<\infty$'' is replaced by ``$\Gdim_{R_\p}M_\p<\infty$'' for every $\p\in\spec R$ with $\depth R_\p<n$; see \cite[Theorem 41]{Masek} and \cite[Theorem 2.12(b)]{DS}, respectively.
In the proof of the above corollary, the local finiteness of the G-dimension of the dual module implied the grade conditions required in Corollaries \ref{EGhighAB} and \ref{Serregrade}, respectively, so we could directly apply these corollaries. It is then natural to ask whether the finiteness of G-dimension without taking duals also implies the desired grade conditions. More precisely, if $\Gdim_{R_\p}M_\p<\infty$ for every $\p\in\spec R$ with $\depth R_\p<n$, does it follow that $\grade_R\Ext^{j+1}_R(\tr M,R)\ge j$ for every $1\le j\le n$, or that $\grade_R\Ext^k_R(\tr M,R)\ge k$ for every $1\le k\le n$?
When $n=1$, it is easy to see that both conditions hold. Indeed, for every associated prime $\p$ of $R$, the assumption on the finiteness of G-dimension and the Auslander--Bridger formula imply that $M_\p$ is a totally reflexive $R_\p$-module, and hence so is $\tr M_\p$. In particular, we have $\Ext^1_{R_\p}(\tr M_\p,R_\p)=\Ext^2_{R_\p}(\tr M_\p,R_\p)=0$. This means that the desired grade conditions are satisfied.
However, when $n\ge2$, these conditions do not necessarily hold. For simplicity, consider the case $n=2$, and let $(R,\m,k)$ be a non-Gorenstein local ring of depth one. Set $M=\tr k$. Since $M$ has projective dimension one, it satisfies the desired condition on the local finiteness of G-dimension. However, since $R$ is not Gorenstein, $\Ext^i_R(\tr M,R)\cong\Ext^i_R(k,R)$ is a nonzero $k$-vector space for every $i\ge1$, and hence has grade exactly one. Thus neither of the desired grade conditions is satisfied.
As discussed above, the local finiteness of G-dimension without taking duals does not imply the grade conditions in general. However, if $M$ is an $(n+1)$-st syzygy or satisfies the modified Serre condition $(\widetilde{S_n})$, then the Auslander--Bridger formula implies local total reflexivity at the relevant prime ideals. Consequently, the grade conditions required in Corollaries \ref{EGhighAB} and \ref{Serregrade} are satisfied, and hence higher torsionfreeness follows.
\end{rmk}

\begin{ex}
Let $M$ be a finitely generated module over a commutative noetherian ring $R$ satisfying the modified Serre condition $(\widetilde{S_1})$. If $M$ is not torsionless, then Corollary \ref{Serregrade} implies that $\Ext^1_R(\tr M,R)$ must have grade zero, that is, this Ext module is not torsion.
For example, let $R=k[[x,y,z]]/(x^2,xy,y^2)$, which is a one-dimensional Cohen--Macaulay local ring. 
Consider the $R$-module $M=R/(x)$. Then $M$ is a maximal Cohen--Macaulay $R$-module, and hence satisfies $(\widetilde{S_1})$. However, $M$ is not a torsionless $R$-module. Indeed, $M^*\cong(0:_R x)=(x,y)$ is annihilated by $y$, while $M$ is not. Hence $M$ cannot be embedded into $M^{**}$, and thus is not torsionless.
In this case, we have
$$
\Ext^1_R(\tr M,R)\cong\Ext^1_R(R/(x),R)\cong(x,y)/(x)\cong(y)\cong R/(x,y).
$$
Note that $\Hom_R(R/(x,y),R)\cong(x,y)\ne0$, so this Ext module is not torsion. In other words, $\grade_R\Ext^1_R(\tr M,R)=0$.
\end{ex}

\begin{ex}
Let $R$ be a commutative noetherian ring and set $\Lambda=\begin{pmatrix}R&R\\0&R\end{pmatrix}$. Let $e_{ij}$ denote the standard matrix units, and set $e_1=e_{11}$ and $e_2=e_{22}$. Let $P_1=e_1\Lambda$ and $P_2=e_2\Lambda$, and define $M=\Cok\iota$, where $\iota:P_2\to P_1$ is the natural injection given by left multiplication by $e_{12}$. Then $\Hom_\Lambda(M,P_1)=0=\Hom_\Lambda(M,P_2)$, so $M$ is a torsion $\Lambda$-module. Hence $\Ext^1_{\Lambda^{\mathrm{op}}}(\tr M,\Lambda)\cong M$, which is isomorphic to $R$ as an $R$-module. Thus $\grade_R\Ext^1_{\Lambda^{\mathrm{op}}}(\tr M,\Lambda)=0$, which is consistent with Corollary \ref{Serregrade}.
\end{ex}

\begin{ex}
Let $k$ be a field, $Q=k[[x,y,z,w]]$, and $I=(xy,zw)$. Consider the ring $R=Q/I^2$ and the $R$-module $M=I/I^2$. Then $R$ is a two-dimensional Cohen--Macaulay local ring. Moreover, setting $S=Q/I$, we have $M\cong S^{\oplus 2}$; see, for example, \cite[Theorem 19.9]{Mat}. Since $S$ is a maximal Cohen--Macaulay $R$-module, so is $M$. In other words, $M$ satisfies the modified Serre condition $(\widetilde{S_2})$. Also, since $M$ is an ideal of $R$, it is clearly a torsionless $R$-module.
Since $(0:_R M)=M$, the image of every element of $M^\ast$ is contained in $M$. Therefore, we have 
$$M^\ast=\Hom_R(M,R)\cong\Hom_R(M,M)\cong\Hom_S(M,M)\cong S^{\oplus 4}.$$
By a similar calculation, we obtain $M^{\ast\ast}\cong S^{\oplus 8}$. In particular, $M$ is not a reflexive $R$-module.
On the other hand, since $S$ is a direct summand of $M^\ast$ as an $R$-module, we have $S\in\add\Syz_2(R)$. Hence $M\in\add\Syz_2(R)$, and thus the Auslander--Bridger exact sequence $0\to M\to M^{\ast\ast}\to\Ext^2_R(\tr M,R)\to0$ splits. Therefore, the above calculations yield $\Ext^2_R(\tr M,R)\cong S^{\oplus(8-2)}=S^{\oplus6}$.
In particular, we conclude that $\grade_R\Ext^2_R(\tr M,R)=0$.
\end{ex}

Finally, we consider the generality of the main result of this section. In this section, we assumed that the two-sided noetherian ring $\Lambda$ is a noetherian $R$-algebra in order to use the depth function. However, we can obtain the same conclusions by using a function satisfying the following properties of depth.
Let $\Lambda$ be a two-sided noetherian ring and $\mathcal{S}$ a nonempty set. Let $\mathcal{M}$ be the set of isomorphism classes of finitely generated right $\Lambda$-modules. Assume that a function $\operatorname{dep}(-,-):\mathcal{S}\times\mathcal{M}\to\ZZ_{\ge0}\cup\{\infty\}$ satisfies the following conditions.
\begin{enumerate}[\rm(1)]
\item For every $\p\in\mathcal{S}$, the function $\operatorname{dep}(\p,-)$ satisfies the Depth Lemma, namely, for any exact sequence $0\to M_{-1}\to M_0\to M_1\to0$ of finitely generated $\Lambda$-modules, one has $\operatorname{dep}(\p,M_0)\ge\min\{\operatorname{dep}(\p,M_{-1}),\operatorname{dep}(\p,M_1)\}$ and $\operatorname{dep}(\p,M_i)\ge\min\{\operatorname{dep}(\p,M_{-i})-i,\operatorname{dep}(\p,M_0)\}$ for $i=\pm1$.
\item $\operatorname{dep}(\p,0)=\infty$ for every $\p\in\mathcal{S}$.
\item If $M\ne0$, then $\ass M:=\{\p\in\mathcal{S}\mid\operatorname{dep}(\p,M)=0\}\ne\emptyset$, namely, there exists $\p\in\mathcal{S}$ such that $\operatorname{dep}(\p,M)=0$.
\end{enumerate}
Then the same arguments as in Lemma \ref{syzsn} and Theorem \ref{SerreAss} apply if depth and associated primes are replaced by the above notions, and the same conclusions hold. Of course, if $\Lambda$ is a noetherian $R$-algebra, these results are recovered by setting $\mathcal{S}=\spec R$ and $\operatorname{dep}(\p,M)=\depth_{R_\p}M_\p$ for every $\p\in\spec R$ and every finitely generated $\Lambda$-module $M$. However, to avoid excessive generality, we stated our results in terms of depth over noetherian algebras.

\section{Triviality of extensions of higher torsionfree modules}
In this section, we study extensions of higher torsionfree modules. The question of when the categories of higher torsionfree modules or higher syzygy modules are closed under extensions has been extensively studied, notably by Auslander--Reiten \cite{AR1,AR2}; see also \cite{DT,GT,Huang99,Hu06,Hu09,HI}. For example, as discussed in the previous section, if the category of $n$-torsionfree modules coincides with that of modules satisfying the modified Serre condition $(\widetilde{S_n})$, then it is closed under extensions. Thus, studying extension-closedness amounts to considering a weaker condition than the equality investigated in the previous section.
The following result of Auslander and Reiten \cite[Theorem 1.6]{AR2} is fundamental.
For a finitely generated $\Lambda$-module $X$, the {\em strong grade} of $X$, denoted by $\sgrade_\Lambda X$, is the infimum of the grades of all submodules of $X$.

\begin{thm}[Auslander--Reiten]\label{ARsgrade}
Let $n\ge0$ be an integer and $C$ an $n$-torsionfree $\Lambda$-module. Then the following are equivalent.
\begin{enumerate}[\rm(1)]
   \item One has the inequality $\sgrade_{\Lambda^{\mathrm{op}}}\Ext^1_\Lambda(C,\Lambda)\ge n$.
   \item For any short exact sequence $0\to A\to B\to C\to 0$ in $\mod\Lambda$ with $A$ $n$-torsionfree, the middle term $B$ is also $n$-torsionfree.
   \item For any short exact sequence $0\to P\to E\to C\to 0$ in $\mod\Lambda$ with $P$ projective, the middle term $E$ is also $n$-torsionfree.
\end{enumerate}
\end{thm}

Based on the above theorem, the extension-closedness of the categories of $n$-torsionfree modules and $n$-syzygy modules can be characterized in terms of grade conditions. In the following theorem, Auslander and Reiten \cite[Theorem 0.1]{AR2} proved that conditions (1), (2), and (4) are equivalent and imply condition (5) for any two-sided noetherian ring $\Lambda$. They also proved that (5) implies (4) if $\Lambda$ is a noetherian algebra over a commutative noetherian ring. Subsequently, Huang and Iyama \cite[Theorem 4.1(3)]{HI} showed that (5) implies (4) without assuming that $\Lambda$ is a noetherian algebra. Moreover, Huang \cite[Theorem 3.3]{Huang99} proved that these conditions are equivalent to (3). Although his result was stated for noetherian algebras, this assumption is in fact unnecessary.

\begin{thm}[Auslander--Reiten, Huang, Huang--Iyama]
Let $\Lambda$ be a two-sided noetherian ring and $n\ge0$ an integer. Then the following conditions are equivalent.
\begin{enumerate}
   \item The category $\Syz_k(\Lambda)$ is extension closed for any $1\le k\le n$.
   \item The category $\add\Syz_k(\Lambda)$ is extension closed for any $1\le k\le n$.
   \item The category $\Trf_k(\Lambda)$ is extension closed for any $1\le k\le n$.
   \item One has $\sgrade_{\Lambda^{\mathrm{op}}}\Ext^{k+1}_\Lambda(M,\Lambda)\ge k$ for any $M\in\mod\Lambda$ and $1\le k\le n$.
   \item One has $\grade_{\Lambda}\Ext^{k}_{\Lambda^{\mathrm{op}}}(N,\Lambda)\ge k$ for any $N\in\mod\Lambda^{\mathrm{op}}$ and $1\le k\le n$.
\end{enumerate}
\end{thm}

Condition (4) is called condition $(c_n)$ by Auslander--Reiten \cite{AR1}, and condition $G_n(1)$ or the quasi $n$-Gorenstein condition by Huang--Iyama \cite{HI}. The opposite version of condition (5) is called condition $(b_n)$ by Auslander--Reiten \cite{AR1}, and condition $g_n(0)$ by Huang--Iyama \cite{HI}. These conditions can also be expressed in terms of the $(l,n)$-condition or the weak $(l,n)$-condition introduced in \cite{IyaSym}; see \cite[5.4]{HI}.

\begin{rmk}\label{exclrmk}
   \begin{enumerate}[\rm(1)]
      \item For any $n\ge0$, the category $\Trf_n(\Lambda)$ is extension closed if and only if it is resolving. Here, a subcategory $\X\subset\mod\Lambda$ is called resolving if it contains $\proj\Lambda$ and is closed under finite direct sums, direct summands, extensions, and syzygies. This was proved by Dey and Takahashi \cite[Proposition 5.1]{DT} under the assumption that $\Lambda$ is commutative. However, their proof does not use commutativity.
      \item When $\Lambda=R$ is commutative and satisfies a Serre condition, remarkable results are known which show that it is not necessary to consider all consecutive integers as in the above theorem. Indeed, by \cite[Theorem B]{GT} and \cite[Proposition 5.1(4)]{DT}, let $R$ be a $d$-dimensional commutative noetherian local ring and $0\le n\le d+1$. If $R$ satisfies $(S_n)$, then the following conditions are equivalent: $\Syz_n(R)$ is extension closed; $\Trf_n(R)$ is extension closed; and $R$ satisfies $(G_{n-1})$. Related results without the assumption that $R$ is local are also known; see \cite{DT,GT,MTT} for details.
   \end{enumerate}
\end{rmk}

Theorem \ref{ARsgrade} says that for an $n$-torsionfree $\Lambda$-module $C$, the inclusion $\proj\Lambda\ast C\subset\Trf_n(\Lambda)$ holds if and only if $C$ satisfies a certain strong grade condition. However, as we shall see below, when $n$ is greater than the little finitistic dimension of $\Lambda$, even the intersection $(\proj\Lambda\ast C)\cap\Trf_n(\Lambda)$ becomes ``trivial''. Here, recall that the (right) {\em little finitistic dimension} of $\Lambda$ is defined to be the supremum of the projective dimensions of finitely generated (right) $\Lambda$-modules of finite projective dimension.

\begin{thm}\label{ARsplit}
Let $\Lambda$ be a two-sided noetherian ring, and assume that $\Lambda$ has finite right little finitistic dimension $t$. Suppose that a $\Lambda$-module $C$ is $t$-torsionfree and that there exists an exact sequence $\sigma:0\to P\to E\to C\to0$ in $\mod\Lambda$ with $P$ projective. If $E$ is $(t+1)$-torsionfree, then $\sigma$ splits.
\end{thm}

To prove the above theorem, we need some preparations. First, we recall the notion of spherical modules introduced by Auslander--Bridger. A finitely generated right $\Lambda$-module $M$ is called {\em $n$-spherical} if $\pd_\Lambda M\le n$ and $\Ext^i_\Lambda(M,\Lambda)=0$ for all $1\le i<n$. The following is precisely the first part of \cite[Theorem 1.1]{Ost} and gives a fundamental property of spherical modules. See also \cite{AB,BHHPT} for other properties.

\begin{lem}\label{spher}
Let $n>0$ be an integer. Then there is an equivalence between the stable category of $n$-spherical $\Lambda$-modules and the category $\grd_n(\Lambda^{\mathrm{op}})$, given by the correspondence $M\mapsto\Ext^n_\Lambda(M,\Lambda)$ for an $n$-spherical $\Lambda$-module $M$, with inverse given by $N\mapsto\tr\syz^{n-1}N$.
\end{lem}

As a direct consequence of the above lemma, we obtain the following corollary, which should be well known.

\begin{cor}\label{findim}
Let $\Lambda$ be a two-sided noetherian ring with finite right little finitistic dimension $t$.
If $E\in\mod\Lambda^{\mathrm{op}}$ satisfies $\grade_{\Lambda^{\mathrm{op}}}E>t$, then $E=0$.
\end{cor}

\begin{proof}
By Lemma \ref{spher}, it is enough to see that every $(t+1)$-spherical $\Lambda$-module is projective.
Let $N$ be a $(t+1)$-spherical $\Lambda$-module. Then $\pd_\Lambda N\le t+1$ by definition, and hence $\pd_\Lambda N\le t$ by the assumption on the finitistic dimension. Since $\Ext^i_\Lambda(N,\Lambda)=0$ for all $1\le i\le t$, the module $N$ is projective.
\end{proof}

Next, we consider morphisms represented by monomorphisms and their higher versions. The former notion was introduced by Auslander--Bridger \cite{AB} and further developed by Kato \cite{Kmorp}. A morphism $f:X\to Y$ in $\mod\Lambda$ is said to be {\em represented by monomorphisms} if there exist projective modules $P,Q$ and morphisms $s,t,u$ such that $\begin{pmatrix}f&s\\t&u\end{pmatrix}:X\oplus P\to Y\oplus Q$ is a monomorphism. A higher version of this notion was introduced in \cite{Omorp}. A morphism $f:X\to Y$ is said to be {\em $n$-representable by monomorphisms} if the map $\Ext^i_{\Lambda^{\mathrm{op}}}(\tr f,\Lambda):\Ext^i_{\Lambda^{\mathrm{op}}}(\tr X,\Lambda)\to\Ext^i_{\Lambda^{\mathrm{op}}}(\tr Y,\Lambda)$ is an isomorphism for all $1\le i<n$ and a monomorphism for $i=n$. This condition was called condition $(\operatorname{T}_n)$ in \cite[Definition 3.1]{Omorp}. In fact, by \cite[Proposition 3.6 and Theorem 3.9]{Kmorp} and \cite[Theorem 3.8]{Omorp}, being $1$-representable by monomorphisms is equivalent to being represented by monomorphisms. Moreover, the following result shows that this higher version is natural.

\begin{lem}\label{embedd}
Let $F:\lmod\Lambda\to\K(\proj\Lambda)$ be the embedding given in \cite{Kkern}, where $\K(\proj\Lambda)$ denotes the homotopy category over $\proj\Lambda$. For an integer $n>0$, a morphism $f:X\to Y$ in $\mod\Lambda$ is $n$-representable by monomorphisms if and only if $H^i(\Cone F(\underline{f}))=0$ for all $-1\le i\le n-2$. Here, $\Cone(-)$ denotes the mapping cone of a morphism of complexes.
\end{lem}

\begin{proof}
The case $n=1$ is precisely \cite[Theorem 3.9]{Kmorp}. For general $n$, recall that the functor $F$ is constructed by taking a projective resolution $P_\bullet$ of $X$ and a projective resolution $Q_\bullet$ of $\tr X$, and splicing $P_\bullet$ and $Q_\bullet^\ast$ together; see \cite[Proposition 2.3]{Kkern}. Thus, by this construction, we obtain an exact sequence
\begin{align*}
0\to H^{-1}(C_\bullet)\to\Ext^1_{\Lambda^{\mathrm{op}}}(\tr X,\Lambda)\to\Ext^1_{\Lambda^{\mathrm{op}}}(\tr Y,\Lambda)\to H^0(C_\bullet)\to\cdots\\
\cdots\to H^{i-2}(C_\bullet)\to\Ext^i_{\Lambda^{\mathrm{op}}}(\tr X,\Lambda)\to\Ext^i_{\Lambda^{\mathrm{op}}}(\tr Y,\Lambda)\to H^{i-1}(C_\bullet)\to\cdots,
\end{align*}
where $C_\bullet=\Cone F(\underline{f})$. The assertion follows from this exact sequence.
\end{proof}

The following lemma plays a key role in the proof of the main result of this section.

\begin{lem}\label{tncor}
Let $0\to K\to X\xrightarrow{f}Y\to0$ be an exact sequence in $\mod\Lambda$ and $n\ge0$. Then
\begin{enumerate}[\rm(1)]
   \item The morphism $f$ is a stable isomorphism, namely, $\underline{f}$ is an isomorphism in $\lmod\Lambda$ if and only if $K$ is projective and $\tr f$ is represented by monomorphisms.
   \item Suppose that the following conditions are satisfied
   \begin{itemize}
      \item The morphism $f$ is $n$-representable by monomorphisms.
      \item The module $K$ is $(n+1)$-torsionfree.
   \end{itemize}
   Then one has the inequality $\grade_{\Lambda^{\mathrm{op}}}\Ker(\Ext^1_\Lambda(f,\Lambda):\Ext^1_\Lambda(Y,\Lambda)\to\Ext^1_\Lambda(X,\Lambda))\ge n$.
\end{enumerate}
\end{lem}

\begin{proof}
The assertions follow from \cite[Lemma 3.5(1)]{Omorp} and \cite[Theorem 3.6]{Omorp}, respectively.
\end{proof}

The following corollary shows that a symmetry under duality holds even when the assumption on the left term in Theorem \ref{ARsplit} is weakened.

\begin{cor}\label{210}
Let $\Lambda$ be a two-sided noetherian ring, and assume $\Lambda$ has finite right little finitistic dimension $t$.
If $0\to A\to B\to C\to0$ is an exact sequence in $\mod\Lambda$ such that $A$ is $(t+2)$-torsionfree, $B$ is $(t+1)$-torsionfree, and $C$ is $t$-torsionfree, then its $\Lambda$-dual $0\to C^\ast\to B^\ast\to A^\ast\to0$ is also exact.
\end{cor}

\begin{proof}
Let $f:B\to C$. Since $\Ext^i_{\Lambda^{\mathrm{op}}}(\tr B,\Lambda)=0$ for all $1\le i\le t+1$ and $\Ext^j_{\Lambda^{\mathrm{op}}}(\tr C,\Lambda)=0$ for all $1\le j\le t$, the morphism $f$ is clearly $(t+1)$-representable by monomorphisms. By the assumption on $A$ and Lemma \ref{tncor}(2), we can see that $$\grade_{\Lambda^{\mathrm{op}}}\Ker(\Ext^1_\Lambda(f,\Lambda):\Ext^1_\Lambda(C,\Lambda)\to\Ext^1_\Lambda(B,\Lambda))\ge t+1.$$
However, by the assumption on the little finitistic dimension of $\Lambda$ and Corollary \ref{findim}, we have $\Ker\Ext^1_\Lambda(f,\Lambda)=0$, that is, $\Ext^1_\Lambda(f,\Lambda)$ is injective. This yields the assertion.
\end{proof}

\begin{proof}[Proof of Theorem \ref{ARsplit}]
For an exact sequence $0\to P\to E\xrightarrow{f}C\to0$, the above corollary implies, in particular, that $\tr f$ is represented by monomorphisms. Hence, Lemma \ref{tncor}(1) applies, and the exact sequence splits.
\end{proof}

Finally, we give some applications of Theorem \ref{ARsplit}. The following corollary is a noncommutative version of \cite[Theorem 1.3(2)]{DT}.

\begin{cor}\label{trfgproj}
Let $\Lambda$ be a two-sided noetherian ring and let $n\ge1$ be an integer.
\begin{enumerate}
\item Assume that $\Lambda$ has finite right little finitistic dimension less than $n$.
If $\Trf_n(\Lambda)$ is extension closed, then every module in $\Trf_n(\Lambda)$ is left Ext-orthogonal to $\Lambda$, that is, semi-Gorenstein-projective in the sense of Ringel--Zhang \cite{RZ1,RZ2}.
\item Assume that both $\Lambda$ and $\Lambda^{\mathrm{op}}$ have finite right little finitistic dimension less than $n$.
If both $\Trf_n(\Lambda)$ and $\Trf_n(\Lambda^{\mathrm{op}})$ are extension closed, then every module in these categories is totally reflexive.
\end{enumerate}
\end{cor}

\begin{proof}
Let $t$ be the right little finitistic dimension of $\Lambda$. We first prove (1).
For each integer $i\ge1$, let $\mathcal{F}_i(\Lambda)$ denote the subcategory of $\mod\Lambda$ consisting of modules $M$ such that $\Ext^j_\Lambda(M,\Lambda)=0$ for all $1\le j\le i$.
We show that $\Trf_n(\Lambda)\subseteq\mathcal{F}_\infty(\Lambda):=\bigcap_{i\ge1}\mathcal{F}_i(\Lambda)$.
Let $M\in\Trf_n(\Lambda)$, and consider an extension $0\to\Lambda\to E\to M\to0$.
Since $\Trf_n(\Lambda)$ is extension closed, we have $E\in\Trf_n(\Lambda)$.
As $n>t$, Theorem \ref{ARsplit} shows that this sequence splits. Hence $\Ext^1_\Lambda(M,\Lambda)=0$.
Moreover, since $M$ is $n$-torsionfree and $\Ext^1_\Lambda(M,\Lambda)=0$, its first syzygy $\syz M$ is $(n+1)$-torsionfree, and hence belongs to $\Trf_n(\Lambda)$.
Applying the same argument to $\syz M$, we obtain $\Ext^1_\Lambda(\syz M,\Lambda)=0$, or equivalently, $\Ext^2_\Lambda(M,\Lambda)=0$.
Repeating this argument, we conclude that $M\in\mathcal{F}_\infty(\Lambda)$.
Thus $\Trf_n(\Lambda)\subseteq\mathcal{F}_\infty(\Lambda)$.

We now prove (2).
By (1), we also have $\Trf_n(\Lambda^{\mathrm{op}})\subseteq\mathcal{F}_\infty(\Lambda^{\mathrm{op}})$.
Let $M\in\Trf_n(\Lambda)$.
Since $M\in\mathcal{F}_\infty(\Lambda)$, we see that $\tr M$ is $\infty$-torsionfree.
In particular, $\tr M\in\Trf_n(\Lambda^{\mathrm{op}})$.
Hence the inclusion established above yields $\Ext^i_{\Lambda^{\mathrm{op}}}(\tr M,\Lambda)=0$ for all $i\ge1$.
Thus $M$ is $\infty$-torsionfree.
Together with $\Ext^i_\Lambda(M,\Lambda)=0$ for all $i\ge1$, this shows that $M$ is totally reflexive.
By symmetry, every module in $\Trf_n(\Lambda^{\mathrm{op}})$ is also totally reflexive.
\end{proof}

We now give an application to commutative algebra. The following statements concern syzygies of the residue field of a commutative noetherian local ring. Assertion (2) is precisely the equivalence of conditions (1), (3), and (4) in \cite[Theorem 1.2]{DT}.

\begin{cor}\label{syzkinter}
Let $(R,\m,k)$ be a commutative noetherian local ring of depth $t$.
\begin{enumerate}[\rm(1)]
   \item For $m\ge t$, the intersection $\proj R\ast \syz^{m}k$ and $\Syz_{t+1}(R)$ is trivial, namely, for an object $M\in\proj R\ast \syz^{m}k$, if $M$ is a $(t+1)$-syzygy module, then $M\cong \syz^m k$ in $\lmod R$.
   \item \textup{(Dey--Takahashi)} The ring $R$ is Gorenstein if and only if $\Trf_{t+1}(R)$ is extension closed, if and only if $\Syz_{n}(R)$ is extension closed for some $n>t$.
\end{enumerate}
\end{cor}

\begin{proof}
Recall that the little finitistic dimension of a commutative noetherian local ring $R$ is equal to its depth $t$. We first prove (1). Assume that there exists an exact sequence $0\to F\to M\to\syz^m k\to0$ with $F\in\proj R$. Then $M$ is locally free on the punctured spectrum and is $(t+1)$-torsionfree by the assumption and Corollary \ref{EGhighAB}. Since $\syz^m k$ is also $t$-torsionfree, Theorem \ref{ARsplit} shows that the exact sequence splits. This proves (1).

Next, we prove (2). If $R$ is Gorenstein, then the other conditions clearly hold. Suppose that $\Trf_{t+1}(R)$ is extension closed. By Corollary \ref{trfgproj}, $\syz^m k$ is totally reflexive for any $m>t$, and hence $R$ is Gorenstein. It remains to show that the extension-closedness of $\Syz_n(R)$ for some $n>t$ implies that $R$ is Gorenstein. In this case, for any exact sequence $0\to F\to M\to\syz^n k\to0$ with $F\in\proj R$, the module $M$ is an $(t+1)$-syzygy module by extension-closedness. Thus, the sequence splits by (1). This implies $\Ext^{n+1}_R(k,R)=0$, and hence $R$ is Gorenstein.
\end{proof}
\section{Submodules of Ext-modules and $n$-reflexive modules}
In this section, we study the $n$-reflexivity of modules.
We first recall its definition.
Let $\Lambda$ be a two-sided noetherian ring.
A module $M\in\mod\Lambda$ is called {\em $n$-reflexive} if the unit map $\underline{\varphi_M^n}:M\to\syz^n\tr\syz^n\tr M$ is an isomorphism in the stable module category $\lmod\Lambda$.
This notion was introduced by Auslander and Bridger \cite{AB} and was later studied in \cite{As,Marc}, but seems to have received relatively little attention.
For small values of $n$, we have the following observations.
When $n=1$, $1$-reflexivity is equivalent to $1$-torsionfreeness.
When $n=2$, Auslander and Bridger gave the following characterization (see \cite[p.~143]{AB}).

\begin{prop}[Auslander--Bridger]\label{ext2proj}
Let $\Lambda$ be a two-sided noetherian ring and $M\in\mod\Lambda$ a torsionless module.
Then $M$ is $2$-reflexive if and only if $\Ext^2_{\Lambda^{\mathrm{op}}}(\tr M,\Lambda)$ is a projective $\Lambda$-module.
\end{prop}

Motivated by the above proposition, we consider the behavior of modules of the form $\Ext^k_{\Lambda^{\mathrm{op}}}(-,\Lambda)$.
Auslander already studied such modules in \cite{Aus69}.
We consider the following subcategories:
\begin{align*}
\E_n^{\mathrm{s}}(\Lambda)&=\{M\in\mod\Lambda\mid M\hookrightarrow\Ext^n_{\Lambda^{\mathrm{op}}}(N,\Lambda)\text{ for some }N\in\mod\Lambda^{\mathrm{op}}\},\\
\E_n(\Lambda)&=\{M\in\mod\Lambda\mid M\cong\Ext^n_{\Lambda^{\mathrm{op}}}(N,\Lambda)\text{ for some }N\in\mod\Lambda^{\mathrm{op}}\}.
\end{align*}

Auslander gave explicit descriptions of $\E_1^{\mathrm{s}}(\Lambda)$ and $\E_1(\Lambda)$  in \cite[Proposition 3.2]{Aus69}.
Higher versions of Auslander's characterizations, namely those for $\E_n^{\mathrm{s}}(\Lambda)$ and $\E_n(\Lambda)$ , are given as follows.
Moreover, modules of high grade provide typical examples of objects in these subcategories.
Here, we recall Auslander's $n$-Gorenstein condition.
The ring $\Lambda$ is called {\em $n$-Gorenstein} in the sense of Auslander if $\sgrade_{\Lambda^{\mathrm{op}}}\Ext^i_\Lambda(M,\Lambda)\ge i$ for every $M\in\mod\Lambda$ and $1\le i\le n$.
This condition is called $(a_n)$ in \cite{AR1} and $G_n(0)$ in \cite{HI}.

\begin{prop}\label{Enexpl}
Let $n>0$ be an integer and $M\in\mod\Lambda$. The following statements hold.
\begin{enumerate}[\rm(1)]
   \item The module $M$ belongs to $\E_n^{\mathrm{s}}(\Lambda)$ if and only if there exists an injection $\iota:M\to X$ in $\mod\Lambda$ such that the induced map $\iota^\ast:X^\ast\to M^\ast$ is zero and $\tr X\in\Syz_{n-1}(\Lambda^{\mathrm{op}})$.
   \item The module $M$ belongs to $\E_n(\Lambda)$ if and only if there exists an injection $\iota:M\to X$ in $\mod\Lambda$ such that the induced map $\iota^\ast:X^\ast\to M^\ast$ is zero, $\tr X\in\Syz_{n-1}(\Lambda^{\mathrm{op}})$, and $\Cok\iota$ is torsionless.
   \item The inclusions $\grd_n(\Lambda)\subset\E_n(\Lambda)\subset\E_n^{\mathrm{s}}(\Lambda)$ hold. Moreover,
   \begin{enumerate}[\rm(i)]
      \item The equalities $\grd_k(\Lambda)=\E_k(\Lambda)$ hold for all $1\le k\le n$ if and only if $\Lambda$ satisfies condition $(b_n)$ in the sense of Auslander--Reiten.
      \item The equalities $\grd_k(\Lambda)=\E_k^{\mathrm{s}}(\Lambda)$ hold for all $1\le k\le n$ if and only if $\Lambda$ satisfies condition $(a_n)$ in the sense of Auslander--Reiten, that is, $\Lambda$ is $n$-Gorenstein.
   \end{enumerate}
\end{enumerate}
\end{prop}

\begin{proof}
We prove (1) and (2). First, we recall the following fact from \cite[Appendix]{AB}: for any $Z\in\mod\Lambda$, consider the Auslander exact sequence
$0\to\Ext^1_{\Lambda^{\mathrm{op}}}(\tr Z,\Lambda)\xrightarrow{j_Z}Z\xrightarrow{\sigma_Z}\D(Z)\to0$,
where $\D(Z)$ is the image of the natural map $Z\to Z^{\ast\ast}$.
Then $\D(Z)\cong\syz\tr\syz\tr Z$ in $\lmod\Lambda$, and the $\Lambda$-dual $\sigma_Z^\ast$ is an isomorphism, or equivalently, $j_Z^\ast$ is zero.
First, we prove the ``only if'' parts of (1) and (2).
Assume that there is an injection $u:M\hookrightarrow\Ext^n_{\Lambda^{\mathrm{op}}}(N,\Lambda)$ with $N\in\mod\Lambda^{\mathrm{op}}$, and set $X=\tr\syz^{n-1}N$.
Then $\tr X\cong\syz^{n-1}N\in\Syz_{n-1}(\Lambda^{\mathrm{op}})$ in $\lmod\Lambda^{\mathrm{op}}$.
Moreover, the composition
$\iota:M\xrightarrow{u}\Ext^n_{\Lambda^{\mathrm{op}}}(N,\Lambda)\cong\Ext^1_{\Lambda^{\mathrm{op}}}(\tr X,\Lambda)\xrightarrow{j_X}X$
is injective, and its $\Lambda$-dual is zero since $j_X^\ast=0$.
If $u$ is an isomorphism, then $\Cok\iota\cong\Cok j_X=\D(X)$ is torsionless by the above fact.
This proves the ``only if'' parts of (1) and (2).
We prove the converse implications.
Suppose that there is an exact sequence $0\to M\xrightarrow{\iota}X\to C\to0$ such that $\iota^\ast=0$ and $\tr X\in\Syz_{n-1}(\Lambda^{\mathrm{op}})$.
Since the dual of the map $X\to C$ is an isomorphism, we obtain the following commutative diagram with exact rows:
$$
\xymatrix@R-1pc@C-1pc{
0\ar[rr]&&M\ar[rr]^{\iota}\ar[dd]&&X\ar[rr]\ar[dd]^-{\varphi_X}&&C\ar[rr]\ar[dd]^-{\varphi_C}&&0\\
&&&&&&&&&\\
0\ar[rr]&&0\ar[rr]&&X^{\ast\ast}\ar[rr]^-{\cong}&&C^{\ast\ast}\ar[rr]&&0.
}
$$
Therefore, writing $\tr X=\syz^{n-1}N$, we see that $M$ is a submodule of $\Ker\varphi_X\cong\Ext^1_{\Lambda^{\mathrm{op}}}(\tr X,\Lambda)\cong\Ext^n_{\Lambda^{\mathrm{op}}}(N,\Lambda)$.
Moreover, if $C$ is torsionless, then $M\cong\Ext^n_{\Lambda^{\mathrm{op}}}(N,\Lambda)$.
This proves the remaining implications of (1) and (2).
The first inclusion in (3) follows directly from Lemma \ref{spher}, and the second is trivial.
Moreover, statements (3)(i) and (3)(ii) follow from the definitions.
\end{proof}

Auslander proved in \cite[Proposition 3.3]{Aus69} that both $\E_1^{\mathrm{s}}(\Lambda)$ and $\E_1(\Lambda)$ are closed under taking quotients.
Clearly, $\E_1^{\mathrm{s}}(\Lambda)$ is closed under taking submodules.
Huang showed in \cite[Example 4.2]{Hu09} that neither $\E_1(\Lambda)$ nor $\E_2(\Lambda)$ is closed under taking submodules.
As we will see below, $\E_2(\Lambda)$ is not closed under taking quotients, even when $\Lambda$ is a commutative artinian ring.

\begin{ex}
Let $(R,\m,k)$ be a non-Gorenstein artinian local ring with $\m^2=0$, and let $e$ be the embedding dimension of $R$.
Then an $R$-module $M$ belongs to $\E_2(R)$ if and only if $M\cong k^{\oplus a(e^2-1)}$ for some integer $a\ge0$.
In particular, $\E_2(R)$ is not closed under taking quotients.
Indeed, applying the $R$-dual to the short exact sequence $0\to\m\to R\to k\to0$, we obtain an exact sequence
$0\to R/\Soc R\to\Hom_R(\m,R)\to\Ext^1_R(k,R)\to0$.
Since $\Soc R=\m\cong k^{\oplus e}$, we have $\Hom_R(\m,R)\cong k^{\oplus e^2}$ and hence $\Ext^1_R(k,R)\cong k^{\oplus(e^2-1)}$.
Now, for any $N\in\mod R$, its first syzygy is a $k$-vector space, say $\syz N\cong k^{\oplus a}$.
Thus $\Ext^2_R(N,R)\cong\Ext^1_R(k,R)^{\oplus a}\cong k^{\oplus a(e^2-1)}$.
Conversely, for any nonzero $x\in\m$, we have $\syz(R/(x))\cong k$, so every module of the form $k^{\oplus a(e^2-1)}$ belongs to $\E_2(R)$.
Since $e\ge2$, we have $e^2-1>1$, and hence $k$ is a quotient of a module in $\E_2(R)$ but does not belong to $\E_2(R)$.
\end{ex}

Auslander showed in \cite[Corollary 3.6]{Aus69} that $\E_1^{\mathrm{s}}(R)$ has the following remarkable property over a commutative noetherian ring $R$.
Combining the following theorem with Proposition \ref{ext2proj}, we obtain the following corollary.

\begin{thm}[Auslander]\label{Ausremarkable}
Let $R$ be a commutative noetherian ring. Then $\E_1^{\mathrm{s}}(R)\cap\proj R=\{0\}$.
In other words, if a projective $R$-module $P$ is a submodule of $\Ext^1_R(N,R)$ for some $N\in\mod R$, then $P=0$.
\end{thm}

\begin{cor}[Auslander--Bridger]\label{AB2ref}
Over a commutative noetherian ring, the notions of $2$-reflexivity and reflexivity coincide.
\end{cor}

It is worth noting that Theorem \ref{Ausremarkable} was established in a more general setting for certain noetherian algebras in \cite[Proposition 3.5]{Aus69}.
Moreover, Auslander proved a stronger result for commutative noetherian rings in \cite[Theorem 4.1]{Aus69}: the same conclusion holds even when the second argument of $\Ext^1$ is an arbitrary module, not necessarily finitely generated.
Motivated by the above corollary, we consider the following questions, which may be regarded as implicit in the work of Auslander and Bridger \cite{AB}.

\begin{ques}\label{ques}
Let $\Lambda$ be a two-sided noetherian ring and $n>0$ an integer.
\begin{enumerate}[\rm(1)]
   \item Is every $n$-reflexive $\Lambda$-module $n$-torsionfree?
   \item Is every $n$-reflexive $\Lambda$-module $(n-1)$-reflexive?
\end{enumerate}
\end{ques}

Question \ref{ques}(1) has an affirmative answer when $n=1$, and Question \ref{ques}(2) has an affirmative answer when $n=2$.
Moreover, if Question \ref{ques}(1) has an affirmative answer for some $n$, then so does Question \ref{ques}(2).
When $n=2$ and $\Lambda$ is commutative, Question \ref{ques}(1) has an affirmative answer by the result of Auslander and Bridger mentioned above.
However, Marczinzik \cite{Marc} recently constructed a finite-dimensional (noncommutative) algebra for which Question \ref{ques}(1) has a negative answer when $n=2$.
To the best of our knowledge, neither question has been settled for $n\ge3$, even when $\Lambda$ is commutative.
The aim of this section is to show that both questions have affirmative answers when $n=3$ and $\Lambda$ is commutative.
To this end, our main tool is the following theorem of Auslander and Reiten, which is none other than \cite[Proposition 3.6]{AR2}.
Recall that for any torsionless $\Lambda$-module $M$, there is an exact sequence
$0\to M\xrightarrow{\varphi_M}M^{\ast\ast}\to\Ext^2_{\Lambda^{\mathrm{op}}}(\tr M,\Lambda)\to0$,
where $\varphi_M$ is the natural map and $\varphi_M^\ast$ is an isomorphism.
Since the double dual functor $(-)^{\ast\ast}$ is isomorphic to $\syz^2\tr\syz^2\tr(-)$, the following theorem can be regarded as a higher version of the classical Auslander--Bridger sequence.

\begin{thm}[Auslander--Reiten]\label{EGARprop}
Let $M \in \mod\Lambda$ and let $n \ge 2$ be an integer.  
If $M$ is $(n{-}1)$-torsionfree, then there exists an exact sequence
$$
0 \to M \xrightarrow{f} \syz^n \tr \syz^n \tr M \oplus Q \to \syz^{n-2} \Ext^n_{\Lambda^{\mathrm{op}}}(\tr M, \Lambda) \to 0
$$
such that $Q$ is projective, the dual $f^\ast$ is surjective, and $\underline{f}$ is isomorphic to the unit morphism $\underline{\varphi^n_M} : M \to \syz^n \tr \syz^n \tr M$ in $\lmod\Lambda$.
\end{thm}

\begin{proof}
Although the original statement does not explicitly state that $f^\ast$ is surjective, we can see from the proof that this is indeed the case.
In fact, the projective module $Q$ appearing in the statement is $L_2^\ast$, as stated in the last line of p.~178.
The modules $L_i$ are defined in the sixth line of the proof, and it follows from their definition that the map $L_2\to M^\ast$ is surjective.
Thus the component $M\to L_2^\ast$ of $f$ has surjective dual, and hence $f^\ast$ is surjective.
\end{proof}

As an application, we obtain a characterization of $n$-reflexive modules, which provides a higher version of Proposition \ref{ext2proj}.

\begin{prop}\label{highernref}
Let $n\ge2$ be an integer and $M\in\mod\Lambda$.
Suppose that $M$ is $(n-1)$-torsionfree.
Then $M$ is $n$-reflexive if and only if $\pd_\Lambda\Ext^n_{\Lambda^{\mathrm{op}}}(\tr M,\Lambda)\le n-2$.
\end{prop}

\begin{proof}
The ``if'' part follows since the short exact sequence in Theorem \ref{EGARprop} splits, and hence $f$ is a stable isomorphism.
Conversely, suppose that $M$ is $n$-reflexive.
Since $f$ is a stable isomorphism, the module $C=\Cok f\cong\syz^{n-2}\Ext^n_{\Lambda^{\mathrm{op}}}(\tr M,\Lambda)$ has projective dimension at most one.
Indeed, setting $W=\syz^n\tr\syz^n\tr M\oplus Q$, we have an exact sequence
$$
\Ext^1_\Lambda(W,-)\xrightarrow{\cong}\Ext^1_\Lambda(M,-)\to\Ext^2_\Lambda(C,-)\to\Ext^2_\Lambda(W,-)\xrightarrow{\cong}\Ext^2_\Lambda(M,-),
$$
which implies that $\Ext^2_\Lambda(C,-)=0$.
Moreover, since $f^\ast$ is surjective, we have an exact sequence
$$
0\to\Ext^1_\Lambda(C,\Lambda)\to\Ext^1_\Lambda(W,\Lambda)\xrightarrow{\cong}\Ext^1_\Lambda(M,\Lambda).
$$
Thus $\Ext^1_\Lambda(C,\Lambda)=0$, which, together with $\pd_\Lambda C\le1$, implies that $C$ is projective.
Consequently, $\pd_\Lambda\Ext^n_{\Lambda^{\mathrm{op}}}(\tr M,\Lambda)\le n-2$.
\end{proof}

Motivated by the above proposition, we consider the following question, which can also be regarded as a natural higher version of Theorem \ref{Ausremarkable} due to Auslander.

\begin{ques}\label{ques2}
Let $n>0$ be an integer and $M\in\mod\Lambda$.
If $M$ has projective dimension at most $n-1$ and $M\in\E_n^{\mathrm{s}}(\Lambda)$, then must $M=0$?
In other words, is zero the only $\Lambda$-module of projective dimension at most $n-1$ that is isomorphic to a submodule of $\Ext^n_{\Lambda^{\mathrm{op}}}(N,\Lambda)$ for some $N\in\mod\Lambda^{\mathrm{op}}$?
\end{ques}

When $n=1$ and $\Lambda$ is commutative, the answer is affirmative by Auslander's theorem, while Marczinzik's example \cite{Marc} shows that the answer is negative in general for noncommutative rings.
Moreover, for general $n$, if Question \ref{ques2} has an affirmative answer, then Proposition \ref{highernref} shows that every $(n-1)$-torsionfree, $n$-reflexive $\Lambda$-module is $n$-torsionfree.
The finiteness of homological dimensions of Ext modules has recently been studied in commutative algebra in \cite{HJM,KimExt}, where the finiteness of the homological dimension of $0$th Ext modules plays an essential role.
It is also important here that we consider projective dimension at most $n-1$ for $n$th Ext modules.
For instance, over a discrete valuation ring $R$ with residue field $k$, the module $\Ext^1_R(k,R)\cong k$ is nonzero and has projective dimension one.

For commutative rings, Question \ref{ques2} has an affirmative answer when $n=2$, as shown below.
The key to the proof is the trace ideal introduced by Auslander and Goldman \cite{AG}.
For $M\in\mod\Lambda$, the {\em trace ideal} of $M$ is defined by $\trace(M)=\sum_{f\in\Hom_\Lambda(M,\Lambda)}\operatorname{Im}f$, which is a two-sided ideal of $\Lambda$.
Note that trace ideals were also studied by Auslander in \cite{Aus69}.

\begin{thm}\label{thmpd1}
Let $R$ be a commutative noetherian ring and $M\in\mod R$.
If $M$ has projective dimension at most one and is isomorphic to a submodule of $\Ext^2_R(N,R)$ for some $N\in\mod R$, then $M=0$.
In other words, Question \ref{ques2} has an affirmative answer when $n=2$ and $\Lambda$ is commutative.
\end{thm}

\begin{proof}
First, we consider the case where $R$ is a one-dimensional Cohen--Macaulay local ring.
By passing to the completion, we may assume that $R$ is complete and hence admits a canonical module $\omega$.
Suppose that $M\ne0$. By assumption, there exists an exact sequence $0\to R^{\oplus a}\xrightarrow{A}R^{\oplus b}\xrightarrow{\pi}M\to0$, where all entries of the matrix $A$ belong to the maximal ideal of $R$.
Applying $\Hom_R(\omega,-)$ to the above exact sequence, we obtain an exact sequence
$$0\to(\omega^\ast)^{\oplus a}\xrightarrow{A}(\omega^\ast)^{\oplus b}\xrightarrow{\pi_\ast}\Hom_R(\omega,M).$$
By \cite[Theorem 2.3]{DKT}, the ideal $\trace(\omega)$ annihilates $\Ext^2_R(N,R)$, and hence it also annihilates $M$.
On the other hand, for each $f\in\Hom_R(\omega,R^{\oplus b})$, the image of each component of $f$ is contained in $\trace(\omega)$.
Thus $\pi_\ast(f)=\pi f$ is the zero map, and hence $\pi_\ast=0$.
It follows from Nakayama's lemma that $(\omega^\ast)^{\oplus b}=0$.
Since $\omega^\ast\ne0$, we must have $b=0$, a contradiction.
Therefore, this completes the proof for the case where $R$ is a one-dimensional Cohen--Macaulay local ring.

Next, we prove the general case.
Suppose that $M\ne0$, and take a minimal prime $\p$ of $M$.
Then $M_\p$ has projective dimension at most one over $R_\p$.
If $M_\p$ is free over $R_\p$, then $M_\p$ is a submodule of $\Ext^2_{R_\p}(N_\p,R_\p)$, and hence $M_\p=0$ by Auslander's theorem (Theorem \ref{Ausremarkable}), a contradiction.
Thus $M_\p$ has projective dimension one. Moreover, by the choice of $\p$, the module $M_\p$ has finite length over $R_\p$.
It follows from the new intersection theorem \cite{PS,Rob87} that $R_\p$ is a one-dimensional Cohen--Macaulay local ring.
Therefore, applying the case already proved, we obtain $M_\p=0$, a contradiction.
This completes the proof of the general case.
\end{proof}

We can also give an affirmative answer to Question \ref{ques}(2) when $n=3$ over commutative rings.
This is achieved by comparing the dimensions of socles of modules.
For a finitely generated $R$-module $M$ over a commutative noetherian local ring $R$ with residue field $k$, we put $\operatorname{s}(M)=\dim_k\Hom_R(k,M)$.
Also, we denote by $\mu(M)$ the minimal number of generators of $M$, that is, $\mu(M)=\dim_k(M\otimes_R k)$.

\begin{lem}\label{3refref}
Every $3$-reflexive module over a commutative noetherian ring is reflexive, and hence $2$-reflexive.
\end{lem}

\begin{proof}
Let $M$ be a $3$-reflexive module over a commutative noetherian ring $R$.
Since $M$ is in particular a second syzygy, Lemma \ref{EGaddsyz} shows that the Auslander exact sequence $0\to M\to M^{\ast\ast}\to\Ext^2_R(\tr M,R)\to0$ splits.
We show that $E=\Ext^2_R(\tr M,R)$ is zero.
Suppose that $E\ne0$, and take $\p\in\ass_R E$.
Since $E$ is a direct summand of $M^{\ast\ast}$, it is torsionless. Hence $\p\in\ass_R R$, so $\depth R_\p=0$.
Moreover, $M_\p$ is $3$-reflexive, and the natural map $M_\p\to(M_\p)^{\ast\ast}$ is a split monomorphism with cokernel $E_\p$. Replacing $R$, $M$, and $E$ by $R_\p$, $M_\p$, and $E_\p$, respectively, we may assume that $(R,\m,k)$ is a local ring with $\depth R=0$ and $\Soc E\ne0$.
Now write $M\cong M'\oplus R^{\oplus r}$, where $M'$ has no nonzero free direct summand. Since $3$-reflexivity is preserved under stable isomorphisms, $M'$ is also $3$-reflexive. Moreover, the natural map $M'\to M'^{\ast\ast}$ is a split monomorphism with cokernel isomorphic to $E$. Thus, replacing $M$ by $M'$, we may further assume that $M$ has no nonzero free direct summand.

Let $p:R^{\oplus a}\to M^\ast$ be a minimal free cover of $M^\ast$, and set $K=\Ker p$. Thus $K$ is a minimal syzygy of $M^\ast$.
Similarly, let $q:R^{\oplus b}\to K^\ast$ be a minimal free cover of $K^\ast$, and set $L=\Ker q$.
Thus we have exact sequences
$$
0\to K\to R^{\oplus a}\xrightarrow{p}M^\ast\to0,
\qquad
0\to L\to R^{\oplus b}\xrightarrow{q}K^\ast\to0.
$$
Note that $K$ has no nonzero free direct summands.
Indeed, such a summand would contradict $(\Soc R)K=0$, which follows from the minimality of $p$.
Similarly, $L$ has no nonzero free direct summands.
Now we have $L\cong\syz(K^\ast)\cong\syz^3\tr\syz (M^\ast)\cong\syz^3\tr\syz^3\tr M$ in $\lmod R$.
Since $M$ is $3$-reflexive, it follows that $M$ and $L$ are isomorphic in $\lmod R$.
However, since neither $M$ nor $L$ has nonzero free direct summands, they are isomorphic in $\mod R$; see \cite[Corollary 1.16]{LW}.
Now, for any finitely generated $R$-module $U$, note that $\Hom_R(k,U^\ast)\cong\Hom_R(k\otimes_R U,R)$, and hence $\operatorname{s}(U^\ast)=\mu(U)\operatorname{s}(R)$.
Thus, applying $\Hom_R(k,-)$ to the first short exact sequence, we obtain $\operatorname{s}(K)=a\operatorname{s}(R)=\mu(M^\ast)\operatorname{s}(R)=\operatorname{s}(M^{\ast\ast})$.
Similarly, the second short exact sequence yields $\operatorname{s}(L)=b\operatorname{s}(R)=\mu(K^\ast)\operatorname{s}(R)=\operatorname{s}(K^{\ast\ast})$.
Since $K$ is torsionless, we have $\operatorname{s}(K)\le\operatorname{s}(K^{\ast\ast})$.
Moreover, since $M\cong L$, we obtain
$$
\operatorname{s}(M)+\operatorname{s}(E)=\operatorname{s}(M^{\ast\ast})=\operatorname{s}(K)\le\operatorname{s}(K^{\ast\ast})=\operatorname{s}(L)=\operatorname{s}(M).
$$
Thus $\operatorname{s}(E)=0$, or equivalently, $\Soc E=0$. This is a contradiction. We are done.
\end{proof}

Combining Proposition \ref{highernref}, Theorem \ref{thmpd1}, and Lemma \ref{3refref}, we obtain the following theorem, which gives an affirmative answer to Question \ref{ques}(1) when $n=3$ and $\Lambda$ is commutative.

\begin{thm}\label{3ref3trf}
Over a commutative noetherian ring, the notions of $3$-reflexivity and $3$-torsionfreeness coincide.
\end{thm}

We summarize some consequences of the results obtained so far in the following corollary.

\begin{cor}\label{directcor}
Let $R$ be a commutative noetherian ring.
\begin{enumerate}[\rm(1)]
   \item If the little finitistic dimension of $R$ is at most one, namely, $\depth R_\m\le 1$ for any maximal ideal $\m$ of $R$, then every $n$-reflexive $R$-module which is $(n-1)$-torsionfree is $n$-torsionfree for all $n\ge1$.
   \item If $M\in\mod R$ has finite projective dimension, then $M\in\E_2^{\mathrm{s}}(R)$ if and only if $M\in\grd_2(R)$.
\end{enumerate}
\end{cor}

\begin{proof}
The assertion (1) directly follows from Proposition \ref{highernref} and Theorem \ref{thmpd1}.
We prove the assertion (2). By Proposition \ref{Enexpl}(3), it is enough to prove the ``only if'' part.
Suppose that $M$ has finite projective dimension and is a submodule of $\Ext^2_R(N,R)$ for some $N\in\mod R$.
For any $\p\in\spec R$ with $\depth R_\p\le1$, the Auslander--Buchsbaum formula yields that $\pd_{R_\p}M_\p\le 1$.
Clearly, $M_\p$ is a submodule of $\Ext^2_{R_\p}(N_\p,R_\p)$.
Thus Theorem \ref{thmpd1} implies that $M_\p=0$.
This means that $M\in\grd_2(R)$.
\end{proof}

We now give an Endo--Goto-type criterion for $n$-reflexivity.
Endo and Goto \cite{EGref} showed that, under certain assumptions, the existence of an isomorphism $M\cong M^{\ast\ast}$ implies that the natural map $M\to M^{\ast\ast}$ is an isomorphism.
We establish a higher version of this result in the stable module category.

\begin{prop}\label{EGtypenref}
Let $\Lambda$ be a two-sided noetherian ring, $n>0$ an integer, and $M\in\mod\Lambda$.
Suppose that one of the following conditions holds:
\begin{enumerate}[\rm(1)]
   \item The stable module category $\lmod\Lambda$ is Krull--Schmidt.
   \item The ring $\Lambda$ is a noetherian $R$-algebra over a commutative noetherian ring $R$.
   \item The ring $\Lambda$ is semilocal, that is, $\Lambda/J(\Lambda)$ is semisimple artinian, where $J(\Lambda)$ denotes the Jacobson radical.
\end{enumerate}
Then $M$ is $n$-reflexive if and only if $M$ and $\syz^n\tr\syz^n\tr M$ are isomorphic in $\lmod\Lambda$.
\end{prop}

To prove the proposition under conditions (2) and (3), we prepare the following lemma.

\begin{lem}\label{directfinite}
Let $\Lambda$ be a two-sided noetherian ring.
Suppose that $\Lambda$ is either a noetherian $R$-algebra over a commutative noetherian ring $R$ or a semilocal ring.
Then, for every $M\in\mod\Lambda$, the ring $\Endl_\Lambda(M)$ is {\em directly finite}, that is, for any $a,b\in\Endl_\Lambda(M)$, the equality $ab=1$ implies $ba=1$.
\end{lem}

\begin{proof}
Put $E=\Endl_{\Lambda}(M)$.
First, suppose that $\Lambda$ is a noetherian $R$-algebra over a commutative noetherian ring $R$.
Suppose that $a,b\in E$ satisfy $ab=1$.
Then left multiplication by $a$ induces a surjective $R$-homomorphism $E\to E$, which is an isomorphism since $E$ is a finitely generated module over the noetherian ring $R$.
Since $a(ba-1)=0$, we obtain $ba=1$.
Next, suppose that $\Lambda$ is semilocal.
Since $M$ is finitely generated, $\End_\Lambda(M)$ is semilocal by \cite[Theorem 3.3]{FH}.
Thus its quotient $E=\Endl_\Lambda(M)$ is also semilocal, and hence directly finite by \cite[Proposition 20.8]{Lam}.
\end{proof}

\begin{proof}[Proof of Proposition \ref{EGtypenref}]
The ``only if'' part is clear.
Conversely, suppose that $M$ and $X=\syz^n\tr\syz^n\tr M$ are isomorphic in $\lmod\Lambda$, and let $\underline{u}:M\to X$ be an isomorphism in $\lmod\Lambda$.
Since $X\in\Syz_n(\Lambda)$, we have $M\in\add\Syz_n(\Lambda)$.
Thus, by Lemma \ref{EGaddsyz}, the unit map $\underline{\varphi_M^n}:M\to X$ is a split monomorphism in $\lmod\Lambda$.
Hence there exists a morphism $\underline{v}:X\to M$ such that $\underline{v\varphi_M^n}=\underline{1_M}$.
Under condition (1), the ring $\Endl_\Lambda(M)$ is directly finite by the Krull--Schmidt property.
Indeed, a split monomorphism from $M$ to itself must be an isomorphism, by the uniqueness of direct sum decompositions.
Under conditions (2) and (3), the same conclusion follows from Lemma \ref{directfinite}.
Now consider the endomorphisms $vu$ and $\underline{u}^{-1}\underline{\varphi_M^n}$ of $M$.
Since $(\underline{vu})(\underline{u}^{-1}\underline{\varphi_M^n})=\underline{1_M}$, direct finiteness yields $(\underline{u}^{-1}\underline{\varphi_M^n})(\underline{vu})=1_M$, and hence $\underline{\varphi_M^n v}=\underline{1_X}$.
Therefore $\underline{\varphi_M^n}$ is an isomorphism in $\lmod\Lambda$, and so $M$ is $n$-reflexive.
\end{proof}
\begin{ac}
The author would like to thank Kaito Kimura and Ryo Takahashi for valuable discussions.
The author is also grateful to Osamu Iyama and Ryu Tomonaga for their useful comments.
ChatGPT was also used to assist with English translation, literature searches, and mathematical discussions.
The author was partly supported by Grant-in-Aid for JSPS Fellows 26KJ0067.
\end{ac}

\end{document}